\documentclass[11pt]{amsart}
\usepackage[T1]{fontenc}
\usepackage{lmodern}
\usepackage{tikz}
\usepackage{xcolor}
\newcommand{\rev}[1]{#1}
\newcommand{\editnov}[1]{#1}
\newcommand{\newrev}[1]{#1}
\usepackage{amsmath,amssymb,amsthm,mathtools}
\usepackage[margin=30mm]{geometry}
\usepackage[colorlinks=true,allcolors=blue]{hyperref}
\newtheorem{theorem}{Theorem}[section]
\newtheorem{corollary}[theorem]{Corollary}
\newtheorem{proposition}[theorem]{Proposition}
\newtheorem{definition}[theorem]{Definition}
\theoremstyle{remark}\newtheorem{remark}[theorem]{Remark}
\newtheorem{example}[theorem]{Example}
\newcommand{\Par}{\mathcal P}
\newcommand{\PP}{\mathcal{PP}}
\newcommand{\Id}{\rev{\mathcal I}}
\newcommand{\Irr}{\rev{\operatorname{Irr}}}
\newcommand{\kk}{\mathbf k}
\newcommand{\NN}{\mathbb N}
\newcommand{\ZZ}{\mathbb Z}
\newcommand{\RR}{\mathbb R}
\DeclareMathOperator{\Spec}{Spec}
\DeclareMathOperator{\Proj}{Proj}
\DeclareMathOperator{\cone}{cone}
\DeclareMathOperator{\Hilb}{Hilb}
\title{Plane Partitions Meet Toric Geometry}
\author{Nancy Abdallah}
\address{Department of Mathematical Sciences, Chalmers University of Technology and the University of Gothenburg, Gothenburg, Sweden}
\email{nancya@chalmers.se}
\author{Mounir Hajli}
\address{Westlake University, Hangzhou 310024, Zhejiang, People's Republic of China}
\email{hajlimounir@gmail.com}
\author{Hussein Mourtada}
\address{Universit\'e Paris Cit\'e, Sorbonne Universit\'e, CNRS, Institut de Math\'ematiques de Jussieu-Paris Rive Gauche, Paris, 75013, France}
\email{hussein.mourtada@imj-prg.fr}
\date{}
\begin{document}
\begin{abstract}
We interpret plane partitions as a basis of global sections of the structure sheaf on a projective limit of normal affine toric varieties. Using this interpretation, we introduce a signed incomparability statistic and prove a product formula for its generating series through Gr\"obner bases and a graph expansion of Hilbert series. The associated projective toric varieties have integral polytopes defined by ordinary partitions. Counting integer points in their dilations counts plane partitions with prescribed row constraints, identifying the corresponding Hilbert functions with Ehrhart polynomials. This connects plane-partition enumeration with toric geometry and Ehrhart theory. \rev{A classical Grassmannian degeneration also yields a geometric proof of MacMahon's formula.} The construction extends recursively to solid and higher-dimensional partitions, independently of an explicit MacMahon-type formula.
\end{abstract}
\maketitle

\section{Introduction}
An (ordinary) integer partition is a finite sequence of positive integers in weakly decreasing order. We write it as $\lambda=(\lambda_1,\ldots,\lambda_r)$ and call $|\lambda|=\sum_i\lambda_i$ its weight. The empty partition $\varnothing$ has weight zero, and $p(n)$ denotes the number of partitions of $n$.
Among the most celebrated identities in partition theory are the Rogers--Ramanujan identities, which equate partitions satisfying difference conditions with partitions satisfying congruence conditions \cite{Andrews}. Their connections with statistical mechanics \cite{Baxter}, probability \cite{Fulman}, and representation theory \cite{FF} illustrate the reach of partition identities across mathematics. More recently, their interpretation through the geometry of arc spaces \cite{BMS,Mourtada} has provided both new proofs and new families of identities. Motivated by this interplay between enumeration and geometry, we develop a toric approach to plane partitions.

Plane partitions generalize ordinary partitions by replacing a decreasing sequence with a two-dimensional array. A \emph{plane partition} is a finitely supported array
$\pi=(\pi_{ij})_{i,j\geq1}$ of nonnegative integers satisfying
$\pi_{ij}\geq\pi_{i+1,j}$ and $\pi_{ij}\geq\pi_{i,j+1}$.
Its volume is $|\pi|=\sum_{i,j}\pi_{ij}$. For instance, the plane partitions of $3$ are \[
\begin{pmatrix}3\end{pmatrix},
\qquad
\begin{pmatrix}2&1\end{pmatrix},
\qquad
\begin{pmatrix}2\\1\end{pmatrix},
\qquad
\begin{pmatrix}1&1&1\end{pmatrix},
\qquad
\begin{pmatrix}1\\1\\1\end{pmatrix},
\qquad
\begin{pmatrix}1&1\\1&0\end{pmatrix}.
\]
Trailing zeros are usually omitted; here we retain them when displaying a rectangular matrix. After removing these zeros, every nonzero row and column is an ordinary partition.

Denote the set of plane partitions by $\PP$.
The classical theory goes back to MacMahon's work on partitions \cite{MacMahon1896,MacMahon}; Stanley's two-part account \cite{StanleyPP} and the treatments in \cite{Andrews,StanleyEC2} describe its subsequent development and its connections with symmetric functions. MacMahon's formula states that
\begin{equation}\label{eq:macmahon}
 \sum_{\pi\in\PP}q^{|\pi|}=\prod_{n\geq1}(1-q^n)^{-n}.
\end{equation}
Plane partitions also connect symmetric functions, probability, and geometry. Okounkov and Reshetikhin \cite{OR} study random plane partitions through the Schur process, while Okounkov, Reshetikhin, and Vafa \cite{ORV} relate crystal melting to toric Calabi--Yau geometry and the topological vertex. Vuleti\'c \cite{Vuletic} gives a two-parameter generalization of MacMahon's formula using Macdonald symmetric functions.

Recent developments include Amanov and Yeliussizov's recovery of MacMahon-type products for higher-dimensional partitions using corner-hook volume \cite{AY}, and Schreier-Aigner's work on symmetry classes of plane and higher-dimensional partitions \cite{SA}. In a related algebraic direction, Graffeo, Monavari, Moschetti, and Ricolfi use socle reduction to study higher-dimensional partition enumeration \cite{GMMR}. These works concern different statistics or restricted classes, rather than replacing MacMahon's ordinary-volume formula for unrestricted plane partitions.

Our approach starts from a simple observation: the rows of a plane partition are ordinary partitions ordered by containment. We use this structure to realize plane partitions as a basis of global regular functions on a projective limit of toric varieties. At each finite stage, diagrams of ordinary partitions define a polytope whose integer points count plane partitions with prescribed row constraints. This construction is recursive: plane partitions provide the corresponding data for solid partitions, and the same procedure continues in higher dimensions. It also leads to a signed identity reflecting the incomparability relations among ordinary partitions.

We first state this identity. It involves a signed incomparability statistic, whose motivation through the gcd--lcm reconstruction problem is explained below.
\begin{definition}\label{def:defect}
Let $\pi\ne\varnothing$ be written on its rectangle of nonzero rows and columns. Let $\mathcal M(\pi)$ consist of matrices with the same column multisets as $\pi$, whose rows are weakly decreasing, pairwise distinct, and arranged in decreasing lexicographic order.
For $A\in\mathcal M(\pi)$, let $\rev{\mathcal E}(A)$ be the set of unordered pairs of incomparable rows, where rows are compared componentwise. Put
\[
 S(\pi)=\sum_{A\in\mathcal M(\pi)}
 \sum_{\substack{D\subseteq\rev{\mathcal E}(A)\\
                  \text{every row occurs in a pair of }D}}(-1)^{|D|},
 \qquad S(\varnothing)=1.
\]
\end{definition}
The lexicographic order counts each row set once and automatically leaves the first column in its original decreasing order. A row comparable with all other rows makes the corresponding inner sum zero.

\begin{example}\label{ex:firststatistic}
Consider the plane partition of volume nine
\[
 \pi=\begin{pmatrix}3&1&1\\2&1&0\\1&0&0\end{pmatrix}.
\]
Its first-column entries are distinct, so the canonical row order is fixed. There are three choices for the zero in column two; the one in column three must then occupy one of the other two rows. Thus $\mathcal M(\pi)$ consists of the following six matrices:
\begin{align*}
 A_1&=\begin{pmatrix}3&0&0\\2&1&1\\1&1&0\end{pmatrix},&
 A_2&=\begin{pmatrix}3&0&0\\2&1&0\\1&1&1\end{pmatrix},&
 A_3&=\begin{pmatrix}3&1&1\\2&0&0\\1&1&0\end{pmatrix},\\
 A_4&=\begin{pmatrix}3&1&0\\2&0&0\\1&1&1\end{pmatrix},&
 A_5&=\begin{pmatrix}3&1&1\\2&1&0\\1&0&0\end{pmatrix},&
 A_6&=\begin{pmatrix}3&1&0\\2&1&1\\1&0&0\end{pmatrix}.
\end{align*}
For $A_1$, the incomparable pairs of row indices are $\{1,2\}$ and $\{1,3\}$. Both must be selected to cover all rows, giving contribution $(-1)^2=1$. For $A_4$, the incomparable pairs are $\{1,3\}$ and $\{2,3\}$, again giving $1$.
In $A_2$, all three pairs are incomparable. Each of the three two-pair selections covers the rows, as does the selection of all three pairs. Its contribution is therefore
\[
 3(-1)^2+(-1)^3=2.
\]
In $A_3$, the first row is comparable with both others; in $A_5$ and $A_6$, the last row is comparable with both others. No selection covers these rows, so these matrices contribute zero. Altogether,
\[
 S(\pi)=1+2+0+1+0+0=4.
\]
In particular, the contribution of a matrix depends on all covering selections of incomparable pairs, not just on the parity of the total number of incomparable pairs.
\end{example}

The motivation for this definition begins with an elementary inverse problem: given positive integers $d,m$, determine the unordered pairs $\{a,b\}$ such that
$\gcd(a,b)=d$ and $\operatorname{lcm}(a,b)=m$.
There are no solutions unless $d\mid m$. Otherwise, writing $a=du$, $b=dv$, this is equivalent to
$uv=m/d$ and $\gcd(u,v)=1$. Each prime-power factor of $m/d$ must be assigned entirely to one of $u,v$.
If $d=m$, the sole solution has $a=b=d$ and hence is not a two-element set if ``pair'' is understood strictly.
Under divisibility, gcd and lcm are meet and join. For two ordinary partitions, their analogues are componentwise minimum and maximum: prescribing both specifies the two entries in each column as an unordered pair. Recovering two decreasing rows is therefore an analogous reconstruction problem. The canonical row order removes the ordering ambiguity. For several rows, we select incomparable pairs covering every row and take their alternating sum.

More explicitly, if $v_\ell(a)$ denotes the exponent of a prime $\ell$ in $a$, then
\[
 v_\ell(\gcd(a,b))=\min\{v_\ell(a),v_\ell(b)\},\qquad
 v_\ell(\operatorname{lcm}(a,b))=\max\{v_\ell(a),v_\ell(b)\}.
\]
Thus both problems recover two vectors from their coordinatewise minimum and maximum. For partitions, the additional condition is that both recovered vectors must be weakly decreasing. If $\alpha\leq\beta$ are prescribed, the two-row plane partition with rows $\beta,\alpha$ records exactly these columnwise minima and maxima. Its admissible rearrangements recover pairs $\lambda,\mu$ with
$\lambda\wedge\mu=\alpha$ and $\lambda\vee\mu=\beta$.
For two distinct rows, a covering set of incomparable pairs exists exactly when the rows are incomparable; it then consists of their unique pair and has cardinality one. Thus an incomparable reconstruction contributes $-1$, whereas a comparable reconstruction contributes zero. The comparable pair is $\{\alpha,\beta\}$ itself, just as $\{d,m\}$ is the comparable arithmetic solution under divisibility. For more rows, one sums over all covering selections of incomparable pairs.

\begin{example}[Two inverse problems]
For $d=2$ and $m=12$, the unordered solutions are $\{2,12\}$ and $\{4,6\}$. The first is comparable under divisibility and the second is not.
Similarly, the partition pairs with meet $(1)$ and join $(2,1)$ are $\{(2,1),(1)\}$ and $\{(2),(1,1)\}$. Both have meet $(1)$ and join $(2,1)$; the former is comparable and the latter is incomparable. They illustrate the same reconstruction problem in the lattice of partitions.
\end{example}

\begin{theorem}\label{thm:main}
With $S(\pi)$ as in Definition~\ref{def:defect}, and $S(\varnothing)=1$, one has the identity of formal power series
\begin{equation}\label{eq:target}
 \sum_{\pi\in\PP} S(\pi)q^{|\pi|}
 =\prod_{n\geq1}(1-q^n)^{p(n)-n}.
\end{equation}
\end{theorem}
The identity also has a graph-theoretic meaning. Consider the graph $G_\infty$ whose vertices are nonempty ordinary partitions and whose edges join incomparable pairs. For instance, $(2)$ and $(1,1)$ are joined, whereas $(2,1)$ and $(1)$ are not. The left-hand side can be read as a signed enumeration of its finite subgraphs without isolated vertices: each subgraph contributes a sign determined by its number of edges and a weight equal to the sum of its vertex weights. The theorem shows that these contributions combine into a product governed solely by the ordinary partition numbers. The statistic $S(\pi)$ records this cancellation more finely, grouping subgraphs according to the plane partition obtained by sorting the columns of their vertex labels. Thus the identity gives global enumerative information about $G_\infty$. Its proof is given in Section~\ref{sec:proof}.

\begin{remark}[Convention on zero rows]\label{rem:zerorows}
The convention $S(\varnothing)=1$ supplies the constant term in \eqref{eq:target}. For nonempty plane partitions we always use the rectangle of nonzero rows and columns specified above: adjoining a zero row would prevent any selection of incomparable pairs from covering all rows. We omit all-zero rows and columns, retaining the trailing zeros within this rectangle when its columns are permuted.
\end{remark}

The geometric motivation comes from the relationship between ordinary partitions and arc spaces, surveyed in \cite{Mourtada}. Bruschek, Mourtada, and Schepers \cite{BMS} interpret the Rogers--Ramanujan identities through Hilbert--Poincar\'e series of arc-space coordinate rings. This interplay between partition identities and commutative algebra is developed further in Afsharijoo's work on even--odd partition identities \cite{Afsharijoo}, in the work of Afsharijoo and Mourtada on partition identities and infinite-dimensional Gr\"obner bases \cite{AM}, and in the new companions to the Andrews--Gordon identities obtained by Afsharijoo, Dousse, Jouhet, and Mourtada \cite{ADJM}. Further developments include the work of Barman, Ghosh, and Singh on the G\"ollnitz--Gordon--Andrews identities \cite{BGS}, and that of Ghosh and Barman on a $J$-generalization of the Rogers--Ramanujan--Gordon identities \cite{GB}. \editnov{Afsharijoo and Mourtada \cite{AMNeighborly} use neighborly partitions and hypergraphs to obtain identities dual to the Andrews--Gordon identities. Afsharijoo, Gonz\'alez P\'erez, and Mourtada \cite{AGPM} relate identities for three-colored partitions to the arc Hilbert--Poincar\'e series of $A_r$ surface singularities.} These results illustrate how algebraic interpretations can both recover classical identities and suggest new families of partition identities.

The starting point is elementary: over a field $\kk$, an ordinary partition $\lambda=(\lambda_1,\ldots,\lambda_r)$ corresponds uniquely to the monomial
\[
 z_{\lambda_1}\cdots z_{\lambda_r}
 \quad\text{in}\quad \kk[z_1,z_2,\ldots],
 \qquad \deg z_i=i.
\]
Its degree is the weight of the partition. Suitable quotients of this polynomial ring admit monomial bases indexed by restricted families of partitions; their Hilbert series then count those families. This is the algebraic viewpoint that we extend to plane partitions; see also \cite{StanleyHilbert,MS,HH} for its broader commutative-algebraic setting.

Since each row of a plane partition is an ordinary partition, we introduce one variable for each nonempty ordinary partition and consider
\[
 B=\kk[x_\lambda:\lambda\text{ is a nonempty ordinary partition}],
 \qquad \deg x_\lambda=|\lambda|.
\]
A plane partition determines the product of the variables corresponding to its nonzero rows. For example,
\[
 \begin{pmatrix}2&1\\1&0\end{pmatrix}
 \quad\longleftrightarrow\quad x_{(2,1)}x_{(1)}.
\]
An arbitrary product of variables need not have rows that can be arranged into a plane partition: the partitions $(2)$ and $(1,1)$ are incomparable componentwise. The natural straightening relation
\[
 x_{(2)}x_{(1,1)}=x_{(2,1)}x_{(1)}
\]
replaces this pair by its componentwise maximum and minimum. More generally, we impose
\[
 x_\lambda x_\mu=x_{\lambda\wedge\mu}x_{\lambda\vee\mu},
\]
where $\wedge$ and $\vee$ denote componentwise minimum and maximum, after adjoining trailing zeros when necessary. The resulting quotient, denoted by $A_\infty$, has a basis consisting of products whose partition indices form a decreasing chain. These basis elements are therefore indexed precisely by plane partitions, with degree equal to volume.

Geometrically, $A_\infty$ is the coordinate ring of a scheme $U_\infty$, so plane partitions index a basis of its global regular functions. We construct $U_\infty$ as a projective limit of finite-dimensional normal affine toric varieties associated with finite distributive lattices of ordinary partitions. We determine their dimensions and the convex cones defining them. These finite constructions are rooted in Hibi's theory of distributive lattices, affine semigroup rings, and algebras with straightening laws \cite{Hibi}, and in the description of the associated toric varieties and their cones by Lakshmibai and Mukherjee \cite{LM}.

This interpretation also supplies the proof strategy for Theorem~\ref{thm:main}. Hibi's Gr\"obner basis theorem replaces the straightening relations by monomial relations indexed by incomparable pairs of ordinary partitions, without changing the Hilbert series. The graph expansion of Mohsen and Mourtada \cite{MM} then expresses its numerator as a signed sum over finite graphs without isolated vertices. Related Hilbert-series formulas for monomial ideals appear in Peeva \cite{Peeva} and Goodarzi \cite{Goodarzi}; the precise form needed here is Lemma~2.4 of Mohsen and Mourtada \cite{MM}. Grouping these contributions by the plane partition obtained by sorting their columns gives exactly the statistic $S(\pi)$. MacMahon's formula evaluates the remaining Hilbert series and yields the stated product. \rev{Section~\ref{sec:rectangles} gives a geometric proof of this evaluation using the classical Grassmannian--Hibi degeneration \cite{GL96,GL97}.}

The same geometric viewpoint turns constrained plane-partition enumeration into lattice-point counting. The defining ideals of the finite affine varieties are homogeneous when every variable has degree one. Taking $\Proj$ of their coordinate rings gives projective toric varieties whose integral polytopes are described explicitly by diagrams of ordinary partitions. Theorem~\ref{thm:polytope} identifies integer points in the $s$-th dilation of such a polytope with plane partitions having exactly $s$ nonzero rows, each belonging to the chosen finite lattice. Volume is the sum of the coordinates of the corresponding point. For example, the ordinary partitions contained in $(2,1)$ define a square: its $s$-th dilation has $(s+1)^2$ integer points, hence there are $(s+1)^2$ plane partitions with exactly $s$ nonzero rows, each contained in $(2,1)$.

This correspondence places the counting problem in Stanley's framework of order polytopes \cite{StanleyPolytopes} and in Ehrhart theory \cite{BR}: the count is a polynomial in the dilation factor, and is also the Hilbert function of the standard graded toric coordinate ring. Partition enumeration and lattice-point enumeration therefore give two ways to compute the same quantities. Brion's formula \cite{Brion} makes this correspondence concrete by expressing the generating polynomial as a sum of contributions from the vertex cones; we illustrate the resulting volume count on the square. The recursive feature is essential to this viewpoint. Replacing ordinary partitions by plane partitions produces polytopes counting solid partitions with prescribed slice constraints; repeating the construction treats partitions in every fixed dimension. At each step, partitions of one dimension supply the geometric data for counting partitions of the next. These correspondences remain available even when an explicit product formula for the ordinary-volume generating function is unknown, as is the case beyond plane partitions.

Section~\ref{sec:LH} recalls distributive lattices, Hibi varieties, and their cones. Section~\ref{sec:PT} applies these constructions to ordinary partitions and constructs the projective limit. Section~\ref{sec:proof} proves the main identity. Section~\ref{sec:geometry} develops the relation between projective toric varieties, their polytopes, and plane-partition counting. \rev{Section~\ref{sec:rectangles} derives MacMahon's formula from a classical degeneration of Grassmannians to Hibi varieties. On an affine chart, this connects the algebra indexed by plane partitions to a polynomial ring. Flatness preserves the number of basis elements of each volume, and the weights of the polynomial variables give the product formula. Corollary~\ref{cor:rectangular-graph} gives a rectangular version of our signed identity, from which Theorem~\ref{thm:main} follows by coefficientwise stabilization. The two toric gradings further refine this identity by the number of rows; the signed numerator then factors through the descent polynomial of the grid (Proposition~\ref{prop:refined-toric}).} Section~\ref{sec:solid} gives the recursive extension to solid and higher-dimensional partitions. The refinement by the number of slices also extends to every dimension: the signed numerator factors through a positive descent polynomial of the finite exponent poset, and these factorizations admit coefficientwise limits (Proposition~\ref{prop:higher-refined} and \eqref{eq:higher-refined-limit}).

\medskip
\noindent\textbf{Acknowledgments.}
The third author is grateful to Pooneh Afsharijoo for numerous discussions on plane partitions around 2017, when he began working on this subject. The authors also thank Marc Chardin and Bernard Teissier for several discussions during the preparation of this work. AI tools were used for translation and assistance with drafting and editing the text.

\section{Distributive lattices and their toric varieties}\label{sec:LH}
\subsection{Distributive lattices and Hibi varieties}
\begingroup
Fix an algebraically closed field $\kk$ of characteristic zero.
A lattice is a partially ordered set in which every pair $\alpha,\beta$ has a greatest lower bound $\alpha\wedge\beta$, called its meet, and a least upper bound $\alpha\vee\beta$, called its join. It is distributive if meet distributes over join and join distributes over meet. Our starting objects are finite distributive lattices; their elements will later be ordinary partitions.

\begin{example}[A lattice of ordinary partitions]
Ordered by diagram containment, the set
\[
 L=\{\varnothing,(1),(2),(1,1),(2,1)\}
\]
is a distributive lattice. Meet and join are componentwise minimum and maximum, after padding partitions with zeros. Its only incomparable pair is $(2),(1,1)$, with meet $(1)$ and join $(2,1)$.
\end{example}

For a finite distributive lattice $L$, write $\widehat0$ for its least element. Its homogeneous Hibi algebra is
\[
 R_L=\kk[x_\alpha:\alpha\in L]/
 (x_\alpha x_\beta-x_{\alpha\wedge\beta}x_{\alpha\vee\beta}:
 \alpha\parallel\beta),
 \qquad \deg x_\alpha=1,
\]
where $\alpha\parallel\beta$ means that $\alpha$ and $\beta$ are incomparable.
Hibi's theorem gives a basis consisting of the standard monomials
$x_{\alpha_1}\cdots x_{\alpha_s}$ with
$\alpha_1\geq\cdots\geq\alpha_s$, including $1$ for $s=0$.
The straightening binomials form a Gr\"obner basis for a suitable reverse lexicographic order, with leading terms $x_\alpha x_\beta$ for incomparable indices; see \cite{Hibi,LM,HHO}.
For the lattice in the preceding example, the unique relation is
\[
 x_{(2)}x_{(1,1)}=x_{(1)}x_{(2,1)}.
\]
The affine variety $\Spec R_L$ and the projective variety
$Y_L=\Proj R_L$ are the associated Hibi varieties.
To describe their toric coordinates and dimensions, we now recover a poset from $L$.

An element $u\ne\widehat0$ is \emph{join-irreducible} if
$u=\alpha\vee\beta$ implies $u=\alpha$ or $u=\beta$.
Write $\Irr(L)$ for the set of these elements, with the order inherited from $L$; the least element is excluded.
An \emph{order ideal} of a poset $P$ is a subset $I$ such that
$u\in I$ and $v\leq u$ imply $v\in I$.
We denote the set of order ideals, including the empty ideal, by $\Id(P)$.
Here ``ideal'' always means order ideal, not an ideal in a ring.

\begin{proposition}[Birkhoff's representation]\label{prop:birkhoff}
For a finite distributive lattice $L$, put $Q=\Irr(L)$. The map
\[
 L\longrightarrow\Id(Q),\qquad
 \alpha\longmapsto I_\alpha=\{u\in\Irr(L):u\leq\alpha\}
\]
is a lattice isomorphism, with inverse $I\mapsto\bigvee_{u\in I}u$ and empty join $\widehat0$. In particular,
$I_{\alpha\wedge\beta}=I_\alpha\cap I_\beta$ and
$I_{\alpha\vee\beta}=I_\alpha\cup I_\beta$.
\end{proposition}
\begin{proof}
Every element of a finite lattice is a join of join-irreducibles below it, by induction. In a distributive lattice, a join-irreducible $u$ is join-prime: if $u\leq\alpha\vee\beta$, then
$u=(u\wedge\alpha)\vee(u\wedge\beta)$ forces $u\leq\alpha$ or $u\leq\beta$.
It follows that $I_{\bigvee_{u\in I}u}=I$ for every order ideal $I\subseteq Q$, and that the two displayed maps are inverse lattice homomorphisms.
\end{proof}

\begin{remark}[Starting from a poset]
Conversely, any finite poset $P$ gives a distributive lattice $\Id(P)$, with intersection as meet and union as join. Its join-irreducibles are the principal ideals $\{v\in P:v\leq u\}$, so $\Irr(\Id(P))\simeq P$.
This is the lattice of order ideals associated with $P$.
It should be distinguished from the sublattice generated by a subset of an already specified lattice, as in the construction of $L_n$ below.
\end{remark}

In the five-element example, $Q=\Irr(L)=\{(1),(2),(1,1)\}$, with $(1)<(2)$ and $(1)<(1,1)$. The partition $(2,1)$ corresponds to all of $Q$, whereas $(2)$ corresponds to the ideal $\{(1),(2)\}$.

Using this identification, we may index the generators by ideals $I\in\Id(Q)$, writing $x_I=x_\alpha$ when $I=I_\alpha$ and $x_\varnothing=x_{\widehat0}$. The Hibi algebra has the monomial realization
\begin{equation}\label{eq:hom}
 R_L\simeq\kk\left[w\prod_{u\in I}z_u:I\in\Id(Q)\right]
 \subseteq\kk[w,z_u:u\in Q].
\end{equation}
It is normal, with $\dim R_L=|\Irr(L)|+1=|Q|+1$, and $Y_L$ has dimension $|Q|$.
We also use the dehomogenized algebra
\begin{equation}\label{eq:dehom}
 A_Q=R_L/(x_{\widehat0}-1)
 \simeq\kk\left[\prod_{u\in I}z_u:I\in\Id(Q)\right].
\end{equation}
Its spectrum $U_Q$ is the affine chart $D_+(x_{\widehat0})$ in $Y_L$.
It is a normal affine toric variety of dimension $|Q|$, with a basis indexed by chains in $L\setminus\{\widehat0\}$, including the empty chain.
Thus $\Spec R_L$, $Y_L$, and $U_Q$ have dimensions $|Q|+1$, $|Q|$, and $|Q|$, respectively.
\par\endgroup

\begin{example}[A chain and a branching poset]
If $Q=\{u<v\}$, its ideals are $\varnothing,\{u\},\{u,v\}$. There are no incomparable ideals, so $R_{\Id(Q)}$ is a polynomial ring in three variables and $U_Q\simeq\mathbb A^2_\kk$.
If $Q=\{u,v,w\}$ with $u<v$, $u<w$, and $v\parallel w$, the only incomparable ideals are $\{u,v\}$ and $\{u,w\}$. Writing their coordinates as $X,Y$, and those of $\{u\}$ and $Q$ as $T,Z$, gives
$U_Q=\Spec\kk[T,X,Y,Z]/(XY-TZ)$. This is the three-dimensional quadric cone.
\end{example}

\subsection{\rev{Join-irreducibles and the associated cones}}
\rev{Let $L$ be a finite distributive lattice and put $Q=\Irr(L)$. The cone of its dehomogenized Hibi variety is described by covers and maximal elements of $Q$, following Lakshmibai and Mukherjee} \cite{LM}. We use the dehomogenized convention \eqref{eq:dehom}, with the least lattice element excluded from $\Irr(L)$.
Let $M_Q=\ZZ^Q$, with basis $(f_u)_{u\in Q}$, and let
$N_Q=\operatorname{Hom}(M_Q,\ZZ)$ have dual basis $(e_u)_{u\in Q}$.
For $I\in\Id(Q)$ put $\chi_I=\sum_{u\in I}f_u$. Set
\[
 C_Q=\cone(\chi_I:I\in\Id(Q))\subseteq(M_Q)_\RR,
 \qquad \sigma_Q=C_Q^\vee\subseteq(N_Q)_\RR.
\]
Write $u\lessdot v$ for a cover relation.
\begin{proposition}\label{prop:cone}
The semigroup and the dual cone are given by
\begin{align*}
 C_Q\cap M_Q&=\{a\in\NN^Q:a_u\geq a_v\text{ whenever }u\leq v\}
             =\sum_{I\in\Id(Q)}\NN\chi_I,\\
 \sigma_Q&=\cone\bigl(e_v:v\text{ maximal in }Q;\quad
                         e_u-e_v:u\lessdot v\bigr).
\end{align*}
Consequently $U_Q=\Spec\kk[\sigma_Q^\vee\cap M_Q]$.
\end{proposition}
\begingroup
\begin{proof}
We first determine the integral points of $C_Q$. Let
$a=\sum_{u\in Q}a_u f_u\in C_Q\cap M_Q$. By the definition of $C_Q$, there are real coefficients $b_I\geq0$ such that
\[
 a=\sum_{I\in\Id(Q)}b_I\chi_I
  =\sum_{I\in\Id(Q)}b_I\sum_{u\in I}f_u
  =\sum_{u\in Q}\left(\sum_{\substack{I\in\Id(Q)\\u\in I}}b_I\right)f_u.
\]
Comparing coefficients in the basis $(f_u)_{u\in Q}$ gives
\[
 a_u=\sum_{\substack{I\in\Id(Q)\\u\in I}}b_I\geq0.
\]
Since $a\in M_Q$, these coefficients are integers. If $u\leq v$, every order ideal containing $v$ also contains $u$. Consequently,
\[
 a_u=\sum_{\substack{I\in\Id(Q)\\u\in I}}b_I
 \geq\sum_{\substack{I\in\Id(Q)\\v\in I}}b_I=a_v.
\]
Thus every integral point of $C_Q$ is nonnegative and order-reversing.

Conversely, let $a=\sum_{u\in Q}a_u f_u$ have nonnegative integral coordinates satisfying $a_u\geq a_v$ whenever $u\leq v$. For each integer $h\geq1$, put
\[
 I_h=\{u\in Q:a_u\geq h\}.
\]
If $v\in I_h$ and $u\leq v$, then $a_u\geq a_v\geq h$, so $u\in I_h$. Hence $I_h\in\Id(Q)$. These ideals are empty for all sufficiently large $h$, and
\begin{align*}
 \sum_{h\geq1}\chi_{I_h}
 &=\sum_{h\geq1}\sum_{u\in I_h}f_u\\
 &=\sum_{u\in Q}\#\{h\in\ZZ_{>0}:h\leq a_u\}\,f_u
 =\sum_{u\in Q}a_u f_u=a.
\end{align*}
This expresses $a$ as a nonnegative integral combination of ideal indicators. Together with the preceding inclusion, it proves both equalities describing $C_Q\cap M_Q$.

We next describe the real cone itself. The same coefficient comparison shows that every point of $C_Q$ is nonnegative and order-reversing. Conversely, suppose $a$ has these properties and let
$0<t_1<\cdots<t_r$ be its distinct positive coordinate values, with $t_0=0$. Each set $I_j=\{u\in Q:a_u\geq t_j\}$ is an order ideal, and
\[
 a=\sum_{j=1}^r(t_j-t_{j-1})\chi_{I_j}.
\]
Indeed, a coordinate equal to $t_k$ receives the telescoping sum
$\sum_{j=1}^k(t_j-t_{j-1})=t_k$. The zero vector corresponds to the empty sum. Therefore
\[
 C_Q=\{a\in\RR^Q:a_u\geq0\text{ for all }u,\quad
                         a_u\geq a_v\text{ whenever }u\leq v\}.
\]

It suffices to retain nonnegativity at maximal elements and inequalities along covers. Every element lies below a maximal element, and every comparable pair can be joined by a saturated chain, because $Q$ is finite. Summing the cover inequalities along such chains recovers all the displayed inequalities. In terms of the dual basis, this gives
\[
 C_Q=\left\{a\in(M_Q)_\RR:
 \begin{array}{ll}
 \langle a,e_v\rangle\geq0&\text{for maximal }v,\\
 \langle a,e_u-e_v\rangle\geq0&\text{for }u\lessdot v
 \end{array}\right\}.
\]
Let $K$ be the cone generated by these vectors $e_v$ and $e_u-e_v$ in $(N_Q)_\RR$. The last description says $C_Q=K^\vee$. Since a finitely generated cone is closed, the bipolar theorem gives
$\sigma_Q=C_Q^\vee=(K^\vee)^\vee=K$, as asserted. All cone combinations here have nonnegative real coefficients.

Finally, set $S=C_Q\cap M_Q$. If $m\in M_Q$ and $km\in S$ for some positive integer $k$, then $m=(km)/k\in C_Q$, hence $m\in S$. Thus $S$ is saturated. The monomial realization of $A_Q$ and the integral decomposition proved above identify
\[
 A_Q=\kk[S]=\kk[\sigma_Q^\vee\cap M_Q],
\]
which gives the claimed description of $U_Q$ and proves normality.
\end{proof}
\par\endgroup
To view these as ordinary Euclidean cones, choose a linear extension
$Q=\{u_1,\ldots,u_d\}$, where $d=|Q|$. The bases identify both
$(M_Q)_\RR$ and $(N_Q)_\RR$ with $\RR^d$, with pairing
$\langle a,y\rangle=\sum_{k=1}^d a_k y_k$.
Here $C_Q$ is the cone of exponent vectors, whereas $\sigma_Q$ is the cone defining the toric variety. In coordinates,
\[
 \sigma_Q=\left\{y\in\RR^d:
       \sum_{u_k\in I}y_k\geq0\text{ for every }I\in\Id(Q)\right\}.
\]
The real cone uses nonnegative \emph{real} combinations of its ray generators; the semigroup of monomials consists of lattice points in its dual. These are distinct objects.

\begin{example}[Cones in $\RR^2$ and $\RR^3$]
For $u<v$, in the order $(u,v)$,
\[
 C_Q=\{(a,b):a\geq b\geq0\},\qquad
 \sigma_Q=\cone\{(0,1),(1,-1)\}
          =\{(y_1,y_2):y_1\geq0,\ y_1+y_2\geq0\}.
\]
The two primitive rays form an integral basis, in agreement with $U_Q\simeq\mathbb A^2_\kk$.
For the branching poset $u<v,u<w$, in the order $(u,v,w)$,
\[
 C_Q=\{(a,b,c):a\geq b\geq0,\ a\geq c\geq0\},\quad
 \sigma_Q=\cone\{(0,1,0),(0,0,1),(1,-1,0),(1,0,-1)\}.
\]
This full-dimensional cone has four extremal rays in $\RR^3$, so it is not simplicial; its variety is the singular quadric cone above.
\end{example}

The principal-ideal indicators form an integral basis after a linear extension of $Q$, so the dimension is $|Q|$.
For comparison, the homogeneous cone in $\RR\oplus(M_Q)_\RR$ is
\[
 \widetilde C_Q=\cone\{(1,\chi_I):I\in\Id(Q)\}
 =\{(h,a):h\geq a_u\geq0,\ a_u\geq a_v\ (u\leq v)\}.
\]
For nonempty $Q$, its dual is generated by $e_0-e_u$ for minimal $u$, by $e_v$ for maximal $v$, and by $e_u-e_v$ for covers. It has dimension $|Q|+1$.

\section{The lattice of ordinary partitions and finite toric approximations}\label{sec:PT}
\subsection{Rectangular join-irreducibles}
Let $\Par$ be the set of ordinary partitions, including $\varnothing$.
We pad partitions by zeros and use diagram containment, not dominance order:
$\lambda\leq\mu$ if $\lambda_i\leq\mu_i$ for every $i$.
Then $\Par$ is a distributive lattice with
\[
 (\lambda\wedge\mu)_i=\min(\lambda_i,\mu_i),\qquad
 (\lambda\vee\mu)_i=\max(\lambda_i,\mu_i).
\]
In particular,
$|\lambda|+|\mu|=|\lambda\wedge\mu|+|\lambda\vee\mu|$.
We write $D(\lambda)=\{(i,j):1\leq j\leq\lambda_i\}$ for its diagram, an order ideal in $\ZZ_{>0}^2$ with product order.

\begin{proposition}\label{prop:rect}
The nonempty join-irreducible partitions are exactly the rectangles
$(b^a)=(b,\ldots,b)$ with $a$ parts. If the distinct positive parts of $\lambda$ are $b_1>\cdots>b_s$ with multiplicities $m_1,\ldots,m_s$, then, setting $r_j=m_1+\cdots+m_j$,
\[
 \lambda=\bigvee_{j=1}^s(b_j^{r_j})
\]
is its unique irredundant join decomposition into join-irreducibles.
\end{proposition}
\begin{proof}
Under the diagram correspondence, rectangles are principal ideals of the grid. An ideal is the union of the principal ideals of its maximal elements, uniquely and irredundantly. For $D(\lambda)$ these maximal elements are $(r_j,b_j)$.
\end{proof}

\begin{example}
For $\lambda=(4,4,2,1)$, the distinct positive parts are $4,2,1$, with multiplicities $2,1,1$, and
$\lambda=(4,4)\vee(2,2,2)\vee(1,1,1,1)$.
For $\lambda=(3,1)$ and $\mu=(2,2)$, their meet is $(2,1)$ and their join is $(3,2)$; these partitions are incomparable although they have equal weight.
\end{example}

\subsection{Principal, rectangular, and weight-generated sublattices}
For $\lambda\in\Par$ let
$L_\lambda=\{\mu\in\Par:\mu\leq\lambda\}$.
Then $L_\lambda\simeq\Id(D(\lambda))$, so
\[
 |\Irr(L_\lambda)|=|\lambda|,\quad
 \dim U_{D(\lambda)}=\dim Y_{L_\lambda}=|\lambda|,\quad
 \dim\Spec R_{L_\lambda}=|\lambda|+1.
\]
In particular, $L_{(b^a)}$ consists of partitions with at most $a$ parts, each at most $b$, and its dehomogenized toric variety has dimension $ab$.

\begin{example}[Two finite sublattices]
The principal sublattice
\[
 L_{(2,1)}=\{\varnothing,(1),(2),(1,1),(2,1)\}
\]
has join-irreducibles $(1),(2),(1,1)$, the branching poset of the preceding examples. The rectangular sublattice
\[
 L_{(2,2)}=\{\varnothing,(1),(2),(1,1),(2,1),(2,2)\}
\]
has join-irreducibles $(1),(2),(1,1),(2,2)$, forming a diamond.
These lattices can be embedded in $\RR^2$ by $\mu\mapsto(\mu_1,\mu_2)$, padding with zeros; meet and join are coordinatewise minimum and maximum. This embedding of the \emph{lattice} must be distinguished from the cone of its toric variety, which lies in $\RR^{|\Irr(L)|}$.
For $L_{(2,2)}$, order the cells as $(11,12,21,22)$. Then
\begin{align*}
 C_{D(2,2)}&=\{a:a_{11}\geq a_{12}\geq a_{22}\geq0,
                       \ a_{11}\geq a_{21}\geq a_{22}\},\\
 \sigma_{D(2,2)}&=\cone\{e_{11}-e_{12},e_{11}-e_{21},
                  e_{12}-e_{22},e_{21}-e_{22},e_{22}\}\subseteq\RR^4.
\end{align*}
Writing
\[T=x_{(1)},\quad X=x_{(2)},\quad Y=x_{(1,1)},\quad Z=x_{(2,1)},\quad W=x_{(2,2)},\]
we obtain
\[U_{D(2,2)}=\Spec\kk[T,X,Y,Z,W]/(XY-TZ),\]
the product of a three-dimensional quadric cone with an affine line.
\end{example}

For $n\geq1$, let $L_n$ be the sublattice \emph{generated} by all partitions of weight at most $n$. The generating set itself need not be a lattice: $(2)\vee(1,1)=(2,1)$ already has weight three.
\begin{proposition}\label{prop:ln}
Put $\rho^{(n)}=(\lfloor n/i\rfloor)_{1\leq i\leq n}$. Then
\[
 L_n=L_{\rho^{(n)}},\qquad Q_n:=\Irr(L_n)\simeq
 \{(i,j)\in\ZZ_{>0}^2:ij\leq n\}.
\]
Writing $d(n)=\sum_{i=1}^n\lfloor n/i\rfloor$, one has
\[
 \dim U_{Q_n}=\dim Y_{L_n}=d(n)=\sum_{m=1}^n\tau(m),
 \qquad \dim\Spec R_{L_n}=d(n)+1,
\]
where $\tau(m)$ counts the positive divisors of $m$. Moreover,
$d(n)=n\log n+O(n)$.
\end{proposition}
\begin{proof}
If $|\mu|\leq n$, then $i\mu_i\leq n$, hence $\mu\leq\rho^{(n)}$.
Conversely, every rectangle $(j^i)$ with $ij\leq n$ is among the generators. Every diagram contained in $D(\rho^{(n)})$ is a union of such rectangles, so belongs to the generated lattice. Counting grid points by their product gives the divisor sum. Finally, summing $\lfloor n/i\rfloor=n/i+O(1)$ gives the asymptotic estimate.
\end{proof}

\subsection{The cones of \texorpdfstring{$U_n$}{U-n} in Euclidean coordinates}
Set $A_n=A_{Q_n}$, $U_n=\Spec A_n$, and $\sigma_n=\sigma_{Q_n}$.
Thus $U_n$ is an affine normal toric variety of dimension $d(n)$, not $n$. Its coordinate ring has one generator $x_\lambda$ for each nonempty $\lambda\leq\rho^{(n)}$, with the straightening relations above and 
$x_\varnothing=1$.
Order the cells of $Q_n$ row by row. This fixes identifications
$(M_{Q_n})_\RR\simeq\RR^{d(n)}$ and $(N_{Q_n})_\RR\simeq\RR^{d(n)}$.
The exponent cone is
\[
 C_{Q_n}=\{a\in\RR_{\geq0}^{d(n)}:
 a_{ij}\geq a_{i+1,j},\ a_{ij}\geq a_{i,j+1}
 \text{ whenever both cells occur}\}.
\]
By Proposition~\ref{prop:cone}, in terms of covers and maximal elements,
\begin{equation}\label{eq:sigman}
 \sigma_n=\cone\left(\{e_{ij}-e_{i+1,j},e_{ij}-e_{i,j+1}\}
                 \cup\{e_{ij}:(i,j)\text{ maximal in }Q_n\}\right),
\end{equation}
where differences are included only when both cells occur. A cell is maximal exactly when $(i+1)j>n$ and $i(j+1)>n$.
Equivalently, this embeds the cone by inequalities as
\[
 \sigma_n=\left\{y\in\RR^{d(n)}:
  \sum_{i=1}^n\sum_{j=1}^{\lambda_i}y_{ij}\geq0
  \text{ for every }\varnothing\ne\lambda\leq\rho^{(n)}\right\}.
\]
These cones are full-dimensional, strongly convex, and rational polyhedral. The listed cover and maximal-element vectors are the primitive extremal rays: the corresponding inequalities of $C_{Q_n}$ are irredundant facets.
Their dual lattice points index the monomial functions on $U_n$.

\begin{example}[$n=1$]
Here $\rho^{(1)}=(1)$, $L_1=\{\varnothing,(1)\}$, and $Q_1=\{(1,1)\}$.
Thus $A_1=\kk[x_{(1)}]$, $U_1=\mathbb A^1_\kk$, and
$C_{Q_1}=\sigma_1=\RR_{\geq0}\subseteq\RR$.
\end{example}

\begin{example}[$n=2$]
Here $\rho^{(2)}=(2,1)$, and $L_2=L_{(2,1)}$ is the five-element lattice listed above. In coordinates $(11,12,21)$,
\begin{align*}
 \sigma_2&=\cone\{(0,1,0),(0,0,1),(1,-1,0),(1,0,-1)\},\\
 \sigma_2&=\{(y_1,y_2,y_3):y_1\geq0,\ y_1+y_2\geq0,
                  \ y_1+y_3\geq0,\ y_1+y_2+y_3\geq0\}.
\end{align*}
Explicitly,
\[
 U_2=\Spec\kk[T,X,Y,Z]/(XY-TZ),\qquad
 (T,X,Y,Z)=(x_{(1)},x_{(2)},x_{(1,1)},x_{(2,1)}).
\]
It is three-dimensional and singular only at the origin, as follows by differentiating $XY-TZ$. The nonsimplicial cone has four rays in $\RR^3$.
\end{example}

\begin{example}[$n=3$]
Now $\rho^{(3)}=(3,1,1)$ and
$Q_3=\{11,12,13,21,31\}$, with order relations
$11<12<13$ and $11<21<31$.
The ten elements of $L_3$ are the empty partition and the nine hooks
$(a,1^b)$ with $1\leq a\leq3$, $0\leq b\leq2$.
In the stated coordinate order,
\[
 C_{Q_3}=\{a:a_{11}\geq a_{12}\geq a_{13}\geq0,
                 \ a_{11}\geq a_{21}\geq a_{31}\geq0\},
\]
\[
 \sigma_3=\cone\{e_{13},e_{31},e_{11}-e_{12},e_{12}-e_{13},
                            e_{11}-e_{21},e_{21}-e_{31}\}\subseteq\RR^5.
\]
For instance $e_{12}-e_{13}=(0,1,-1,0,0)$.
Arrange the nine generators as $H=(H_{ab})_{1\leq a\leq3,\,0\leq b\leq2}$, with $H_{ab}=x_{(a,1^b)}$. Meet and join take the minima and maxima of the two indices, so
\[
 A_3=\kk[H_{ab}]/I_2(H),\qquad
 U_3=\{H\in\mathbb A^9_\kk:\operatorname{rank}H\leq1\},
\]
where $I_2(H)$ is generated by all $2\times2$ minors.
This is the five-dimensional affine cone over the Segre embedding of
$\mathbb P^2_\kk\times\mathbb P^2_\kk$. It is smooth away from its vertex and singular at the vertex; its toric cone has six rays in dimension five.
\end{example}

\begin{example}[$n=4$]
Here $\rho^{(4)}=(4,2,1,1)$ and the row-by-row coordinate order is
\[(11,12,13,14,21,22,31,41).\]
Thus $d(4)=8$ and
\begin{align*}
 C_{Q_4}=\{a:\;&a_{11}\geq a_{12}\geq a_{13}\geq a_{14}\geq0,\\
              &a_{11}\geq a_{21}\geq a_{31}\geq a_{41}\geq0,\\
              &a_{12}\geq a_{22}\geq0,\quad a_{21}\geq a_{22}\}.
\end{align*}
The maximal cells are $14,22,41$, and
\begin{align*}
 \sigma_4=\cone\{&e_{14},e_{22},e_{41},
 e_{11}-e_{12},e_{12}-e_{13},e_{13}-e_{14},\\
 &e_{11}-e_{21},e_{21}-e_{31},e_{31}-e_{41},
 e_{12}-e_{22},e_{21}-e_{22}\}\subseteq\RR^8.
\end{align*}
There are eleven extremal rays. The lattice $L_4$ consists of the partitions with at most four parts satisfying
$\lambda_1\leq4$, $\lambda_2\leq2$, $\lambda_3,\lambda_4\leq1$.
There are $26$ such partitions including the empty one, so $U_4$ is the eight-dimensional affine toric variety in $\mathbb A^{25}_\kk$ defined by their straightening binomials. Its torus-fixed point is singular because $\sigma_4$ is not simplicial. For example, it has the relation
$x_{(3)}x_{(2,1)}=x_{(2)}x_{(3,1)}$.
\end{example}

\subsection{The projective limit and functions indexed by plane partitions}
With $A_n$ and $U_n$ as above, since $Q_n$ is an order ideal in $Q_{n+1}$, its ideals remain ideals in $Q_{n+1}$. The monomial realizations therefore give injective homomorphisms
$A_n\hookrightarrow A_{n+1}$ and compatible dominant toric morphisms
$U_{n+1}\to U_n$.
Define
\begin{equation}\label{eq:limit}
 A_\infty=\varinjlim_n A_n,\qquad
 U_\infty=\Spec A_\infty=\varprojlim_n U_n.
\end{equation}
Here ``infinite toric variety'' refers precisely to this affine inverse-limit scheme; it is not a variety of finite type.
Equivalently,
\[
 A_\infty\simeq
 \kk[x_\lambda:\lambda\in\Par]\big/
 \bigl(x_\varnothing-1,\
 x_\lambda x_\mu-x_{\lambda\wedge\mu}x_{\lambda\vee\mu}
 :\lambda\parallel\mu\bigr).
\]
Indeed, every finite collection of partitions belongs to some $L_n$, and the finite straightening presentations are compatible.

\begin{proposition}\label{prop:basis}
Plane partitions index a $\kk$-basis of $A_\infty$. Giving $x_\lambda$ weight $|\lambda|$ defines a grading with finite-dimensional graded pieces, and
\[
 \Hilb(A_\infty;q)=\sum_{\pi\in\PP}q^{|\pi|}
                   =\prod_{n\geq1}(1-q^n)^{-n}.
\]
\end{proposition}
\begin{proof}
The nonzero rows of a plane partition form a chain
$\lambda^{(1)}\geq\cdots\geq\lambda^{(r)}$ of nonempty ordinary partitions.
Associate to it the standard monomial
$x_{\lambda^{(1)}}\cdots x_{\lambda^{(r)}}$; associate $1$ to the empty plane partition. Conversely, such a chain is exactly a plane partition. Every finite chain occurs in some $L_n$, and the standard monomial bases are preserved by the inclusions. The valuation identity for weight makes the relations homogeneous, and the degree is the volume. There are finitely many plane partitions of any fixed volume. The last equality is MacMahon's formula, \rev{proved geometrically in Theorem~\ref{thm:grass-macmahon}.}
\end{proof}
In the monomial model, the function associated with $\pi$ is
$\prod_{i,j}z_{ij}^{h_{ij}}$, where
$h_{ij}=\#\{a:\lambda^{(a)}_i\geq j\}$. This explicitly identifies the row-chain description with an order-reversing exponent array on the grid.

\begin{example}[A plane partition as a function]
For $\pi=\left(\begin{smallmatrix}2&1\\1&0\end{smallmatrix}\right)$, the row chain is $(2,1)\geq(1)$. The associated function on $U_2$ and $U_\infty$ is
$x_{(2,1)}x_{(1)}=z_{11}^2z_{12}z_{21}$, of weight four. Its exponent array is $(h_{11},h_{12},h_{21})=(2,1,1)\in C_{Q_2}\cap\ZZ^3$.
\end{example}

Finally, let $B=\kk[x_\lambda:\lambda\in\Par\setminus\{\varnothing\}]$, graded by weight. Since it has $p(n)$ variables of weight $n$,
\[
 \Hilb(B;q)=\prod_{n\geq1}(1-q^n)^{-p(n)},\qquad
 \frac{\Hilb(A_\infty;q)}{\Hilb(B;q)}
       =\prod_{n\geq1}(1-q^n)^{p(n)-n}.
\]
All products are interpreted coefficientwise in $\ZZ[[q]]$. This explains the product sought in \eqref{eq:target} as a Hilbert-series quotient. A combinatorial interpretation of this quotient by the signed defect sum $S$ is established in Section~\ref{sec:proof}.

\section{Gr\"obner degeneration and the graph expansion}\label{sec:proof}
We prove \eqref{eq:target} by identifying $S(\pi)$ with the contribution of all subgraphs whose vertex sets yield $\pi$ after column sorting.

\subsection{The incomparability graph and the initial ideal}
\begin{figure}[htbp]
\centering
\begin{tikzpicture}[x=4cm,y=3.5cm,every node/.style={font=\normalsize},vertex/.style={draw=black!65,rounded corners=2pt,inner sep=3pt,line width=.5pt}]
\draw[black!40,line width=.5pt] (0.70859,-0.59769) -- (0.34353,-1.00000);
\draw[black!40,line width=.5pt] (0.70859,-0.59769) -- (-0.48211,0.00982);
\draw[black!40,line width=.5pt] (0.70859,-0.59769) -- (0.28682,0.14182);
\draw[black!40,line width=.5pt] (0.34353,-1.00000) -- (-0.02869,-0.31359);
\draw[black!40,line width=.5pt] (0.34353,-1.00000) -- (-0.10709,0.14186);
\draw[black!40,line width=.5pt] (-0.02869,-0.31359) -- (-0.58918,-0.45917);
\draw[black!40,line width=.5pt] (-0.02869,-0.31359) -- (-0.48211,0.00982);
\draw[black!40,line width=.5pt] (-0.02869,-0.31359) -- (-0.28588,0.52348);
\draw[black!40,line width=.5pt] (-0.02869,-0.31359) -- (0.30616,0.66768);
\draw[black!40,line width=.5pt] (-0.02869,-0.31359) -- (0.28682,0.14182);
\draw[black!40,line width=.5pt] (-0.58918,-0.45917) -- (-0.48211,0.00982);
\draw[black!40,line width=.5pt] (-0.58918,-0.45917) -- (-0.10709,0.14186);
\draw[black!40,line width=.5pt] (-0.58918,-0.45917) -- (0.28682,0.14182);
\draw[black!40,line width=.5pt] (-0.48211,0.00982) -- (-0.10709,0.14186);
\draw[black!40,line width=.5pt] (-0.48211,0.00982) -- (-0.15214,0.88581);
\draw[black!40,line width=.5pt] (-0.48211,0.00982) -- (-0.28588,0.52348);
\draw[black!40,line width=.5pt] (-0.10709,0.14186) -- (-0.15214,0.88581);
\draw[black!40,line width=.5pt] (-0.10709,0.14186) -- (-0.28588,0.52348);
\draw[black!40,line width=.5pt] (-0.10709,0.14186) -- (0.30616,0.66768);
\draw[black!40,line width=.5pt] (-0.10709,0.14186) -- (0.28682,0.14182);
\draw[black!40,line width=.5pt] (-0.15214,0.88581) -- (-0.28588,0.52348);
\draw[black!40,line width=.5pt] (-0.15214,0.88581) -- (0.30616,0.66768);
\draw[black!40,line width=.5pt] (-0.15214,0.88581) -- (0.28682,0.14182);
\draw[black!40,line width=.5pt] (-0.28588,0.52348) -- (0.30616,0.66768);
\draw[black!40,line width=.5pt] (-0.28588,0.52348) -- (0.28682,0.14182);
\draw[black!40,line width=.5pt] (0.30616,0.66768) -- (0.28682,0.14182);
\node[vertex,fill=black!5] (v0) at (-1.15000,-1.10000) {$(1)$};
\node[vertex,fill=blue!10] (v1) at (0.70859,-0.59769) {$(2)$};
\node[vertex,fill=blue!10] (v2) at (0.34353,-1.00000) {$(1,1)$};
\node[vertex,fill=green!10] (v3) at (-0.02869,-0.31359) {$(3)$};
\node[vertex,fill=green!10] (v4) at (-0.58918,-0.45917) {$(2,1)$};
\node[vertex,fill=green!10] (v5) at (-0.48211,0.00982) {$(1,1,1)$};
\node[vertex,fill=orange!15] (v6) at (-0.10709,0.14186) {$(4)$};
\node[vertex,fill=orange!15] (v7) at (-0.15214,0.88581) {$(3,1)$};
\node[vertex,fill=orange!15] (v8) at (-0.28588,0.52348) {$(2,2)$};
\node[vertex,fill=orange!15] (v9) at (0.30616,0.66768) {$(2,1,1)$};
\node[vertex,fill=orange!15] (v10) at (0.28682,0.14182) {$(1,1,1,1)$};
\node[anchor=west,font=\small] at (-.94,-1.1) {isolated vertex};
\node[vertex,fill=blue!10,font=\small] at (-0.85,-1.45) {Weight 2};
\node[vertex,fill=green!10,font=\small] at (-0.05,-1.45) {Weight 3};
\node[vertex,fill=orange!15,font=\small] at (0.75,-1.45) {Weight 4};
\end{tikzpicture}
\caption{The induced subgraph of $G_\infty$ on the eleven nonempty ordinary partitions of weight at most four, with all twenty-six edges. Two vertices are adjacent precisely when their partitions are incomparable componentwise, after adjoining trailing zeros. Colors indicate partition weight. The vertex $(1)$ is isolated because it is contained in every nonempty partition. Crossings of edges do not represent vertices.}
\label{fig:incomparability}
\end{figure}
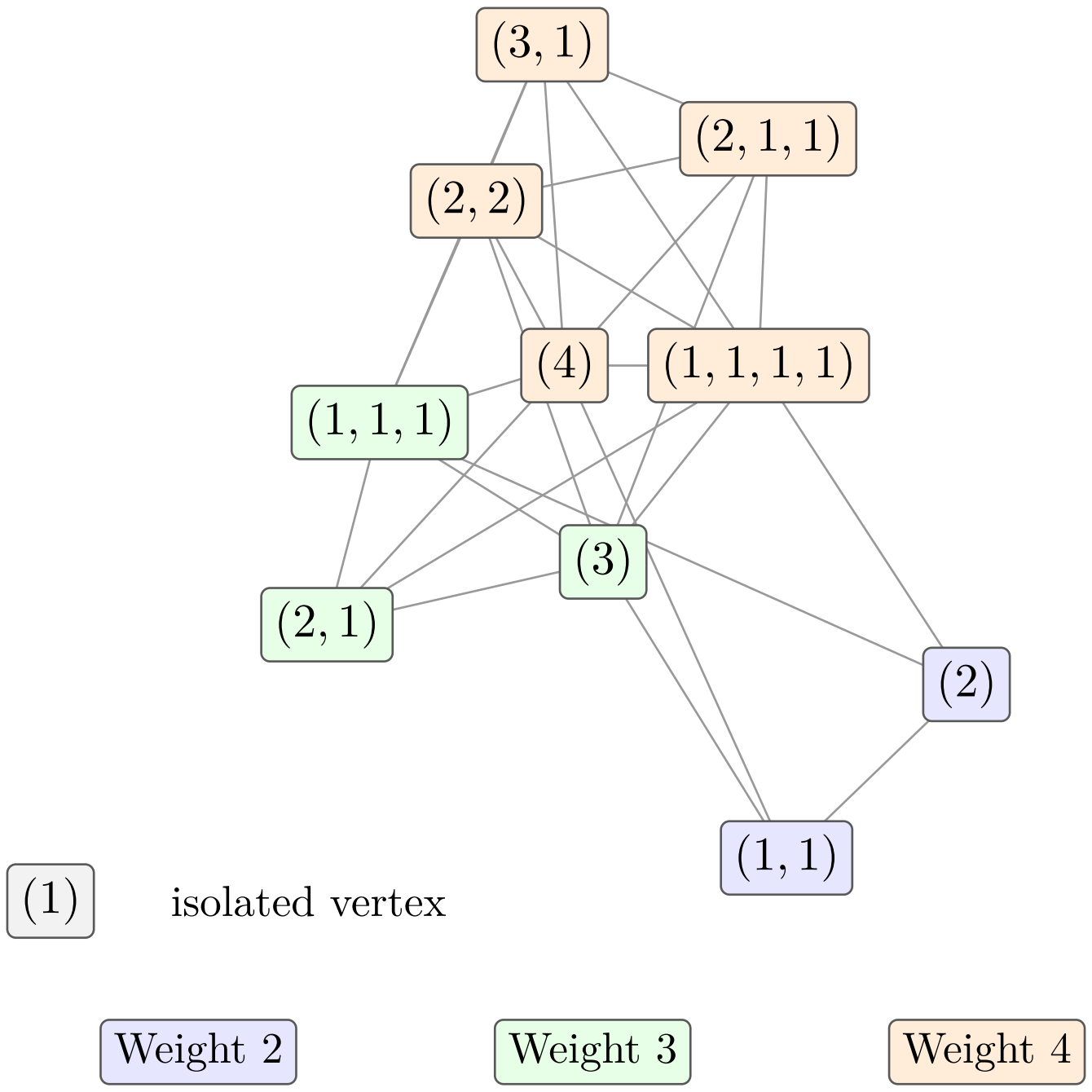

Let $B=\kk[x_\lambda:\lambda\in\Par\setminus\{\varnothing\}]$ with $\deg x_\lambda=|\lambda|$, and let $\mathfrak a$ be the kernel of $B\to A_\infty$.
Define the simple graph $G_\infty$ by
\[
 V(G_\infty)=\Par\setminus\{\varnothing\},\qquad
 E(G_\infty)=\{\{\lambda,\mu\}:\lambda\parallel\mu\}.
\]
Each vertex $\lambda$ has weight $|\lambda|$. Following the standard terminology for monomial ideals \cite{HH,MS}, its edge ideal is
\[
 I(G_\infty)=(x_\lambda x_\mu:\lambda\parallel\mu)\subseteq B.
\]
There are no loops or repeated vertices in this graph. Figure~\ref{fig:incomparability} shows its induced subgraph on nonempty partitions of weight at most four.

\begin{proposition}\label{prop:initial}
There is a weight-compatible monomial order $\prec$ for which
\[
 \operatorname{in}_\prec(\mathfrak a)=I(G_\infty),\qquad
 \Hilb(A_\infty;q)=\Hilb(B/I(G_\infty);q).
\]
\end{proposition}
\begin{proof}
Enumerate the variables in increasing partition weight, resolving ties arbitrarily. This extends the containment order, since strict containment strictly increases weight. Compare monomials first by total weight. For equal weights, use reverse lexicographic comparison: at the largest variable where exponents differ, the monomial with smaller exponent is larger. Each bounded-weight set of monomials is finite, so this defines a well-founded multiplicative order even with infinitely many variables.

For incomparable $\lambda,\mu$, the join is strictly larger than both, and the meet is nonempty. The binomial
\[
 x_\lambda x_\mu-x_{\lambda\wedge\mu}x_{\lambda\vee\mu}
\]
is weight-homogeneous and has initial monomial $x_\lambda x_\mu$.
Hibi's straightening/Gr\"obner basis theorem gives the assertion on every finite $L_n$ (see \cite{HHO}). The relations involving $\varnothing$ are trivial; the nontrivial relations do not involve $x_\varnothing$, since any two nonempty partitions contain $(1)$. Thus the same assertion holds in the dehomogenized presentation.
Every polynomial and every pair of straightening relations lie in some finite $L_n$. The finite Gr\"obner bases consequently give the asserted initial ideal in $B$. Standard monomials preserve weight, proving equality of the graded dimensions, and hence of Hilbert series.
\end{proof}

\begin{example}
The first edge is $\{(2),(1,1)\}$, of total weight four. Its relation is
$x_{(2)}x_{(1,1)}-x_{(1)}x_{(2,1)}$, whose initial term is $x_{(2)}x_{(1,1)}$.
The three partitions $(3),(2,1),(1,1,1)$ form a triangle in $G_\infty$.
\end{example}

\subsection{The subgraph series of Mohsen--Mourtada}
For a simple graph $G$, let $\mathcal H(G)$ consist of its finite subgraphs without isolated vertices, including the empty graph. A subgraph here need not be induced. With independent formal vertex variables $\mathbf u=(u_v)_{v\in V(G)}$, define
\[
 \mathcal S_G(\mathbf u,z)=
 \sum_{H\in\mathcal H(G)}z^{|E(H)|}\prod_{v\in V(H)}u_v.
\]
Related formulas for Hilbert series of monomial ideals occur in the general theory; see Peeva \cite{Peeva} and Goodarzi \cite{Goodarzi}. We use the precise graph-theoretic formulation of Lemma~2.4 of Mohsen--Mourtada \cite{MM}, which expresses the multigraded Hilbert series of the edge-ideal quotient as
\[
 \frac{\mathcal S_G(\mathbf u,-1)}{\prod_{v\in V(G)}(1-u_v)}.
\]
Applying $u_\lambda=q^{|\lambda|}$ and Proposition~\ref{prop:initial} gives
\begin{equation}\label{eq:graphnumerator}
 \sum_{H\in\mathcal H(G_\infty)}
 (-1)^{|E(H)|}q^{\sum_{\lambda\in V(H)}|\lambda|}
 =\prod_{n\geq1}(1-q^n)^{p(n)-n}.
\end{equation}
This specialization is well-defined: only finitely many vertices and edge subsets can contribute to any fixed weight. Thus the finite-graph formula passes coefficientwise to $G_\infty$.

It is possible to group this sum by induced subgraphs, but their weights must be retained. For a finite vertex set $F$, let $G_\infty[F]$ denote the induced graph, and define
\begin{equation}\label{eq:gamma}
 \gamma(G_\infty[F])=
 \sum_{\substack{D\subseteq E(G_\infty[F])\\
                 \text{every vertex of }F\text{ is incident to an edge of }D}}
 (-1)^{|D|}.
\end{equation}
Set $\gamma(\varnothing)=1$. This coefficient is zero if $F$ has an isolated vertex, and \eqref{eq:graphnumerator} can be rewritten as
\[
 \sum_{\substack{F\subseteq V(G_\infty)\\F\text{ finite}}}
 \gamma(G_\infty[F])q^{\sum_{\lambda\in F}|\lambda|}
 =\prod_{n\geq1}(1-q^n)^{p(n)-n}.
\]
The condition on $D$ says that it is an edge cover of the vertex set; it does not require that all induced edges be selected.

\begin{example}[Why induced graphs alone do not give the sign]
For a single edge, $\gamma=-1$. For a path on three vertices, both edges must be selected, so $\gamma=1$. For a triangle, there are three spanning two-edge paths and the full triangle, giving
$\gamma=3(-1)^2+(-1)^3=2$.
Thus $\gamma(K)$ is not in general $(-1)^{|E(K)|}$, even when $K$ has no isolated vertices.
\end{example}

\subsection{Grouping the graph expansion by plane partitions}
Let $F$ be a finite set of nonempty ordinary partitions. Arrange its elements as rows and sort each column in weakly decreasing order. Denote the resulting plane partition by $\pi(F)$.
It is indeed a plane partition: if two columns satisfy the entrywise inequality before sorting, their respective $k$th largest entries satisfy the same inequality after sorting. Its volume is $\sum_{\lambda\in F}|\lambda|$.

Recall that $\mathcal M(\pi)$ consists of the admissible matrices with distinct rows in decreasing lexicographic order. There is a bijection
\[
 \{F\subseteq V(G_\infty):F\text{ finite},\ \pi(F)=\pi\}
 \longleftrightarrow\mathcal M(\pi).
\]
Indeed, sorting the rows of $F$ lexicographically selects exactly one admissible matrix, and taking its row set is the inverse. The distinct-row restriction and the ordering within equal-first-entry blocks are both necessary.
For $A\in\mathcal M(\pi)$ let $K_A$ be the graph on its rows, with an edge for each incomparable pair. Then $K_A=G_\infty[F]$ for the corresponding row set $F$.

By Definition~\ref{def:defect}, the selections of incomparable pairs covering the rows are exactly the edge subsets of $K_A$ covering its vertices. Consequently
\[
 S(\pi)=\sum_{A\in\mathcal M(\pi)}\gamma(K_A)
 \qquad(\pi\ne\varnothing).
\]
This equality translates the matrix definition directly into the coefficient supplied by the subgraph expansion.

\begin{theorem}\label{thm:graphidentity}
The signed global incomparability sums satisfy
\[
 \sum_{\pi\in\PP}S(\pi)q^{|\pi|}
 =\prod_{n\geq1}(1-q^n)^{p(n)-n}.
\]
\end{theorem}
\begin{proof}
Group the finite vertex sets in \eqref{eq:graphnumerator} first by their induced graphs, with coefficient \eqref{eq:gamma}, and then by $\pi(F)$. The bijection above identifies each fiber sum with $S(\pi)$. Volumes are preserved, and coefficientwise finiteness justifies both regroupings.
\end{proof}

\begin{example}[The coefficient in weight six]
The five edges of total vertex weight six are
\[
 \{(2),(1,1,1,1)\},\quad\{(1,1),(4)\},\quad
 \{(3),(2,1)\},\quad\{(3),(1,1,1)\},\quad
 \{(2,1),(1,1,1)\}.
\]
No graph with three distinct vertices and without isolated vertices has weight at most six: $(1)$ is comparable with every partition, and three distinct remaining vertices have total weight at least $2+2+3=7$. Each listed edge contributes $-1$, so
\[
 \sum_{|\pi|=6}S(\pi)=-5
 =[q^6]\prod_{n\geq1}(1-q^n)^{p(n)-n}.
\]
Matrices with repeated rows do not belong to $\mathcal M(\pi)$ and are excluded by definition.
\end{example}

\section{Toric weights, Riemann--Roch, and plane-partition counting}\label{sec:geometry}
We use the standard projective varieties associated with $U_n$ to relate lattice points, sections of line bundles, and plane partitions. Throughout this section, every Proj is taken with respect to the standard grading.

\subsection{Two gradings and sections on the standard Proj}
Every nonempty ideal of $Q_n$ contains its minimum $(1,1)$. Consequently the standard grading $\deg_{\mathrm{std}}x_\lambda=1$ is well-defined on $A_n$, and its value on an exponent vector $a\in C_{Q_n}\cap M_{Q_n}$ is
\[
 \eta_n(a)=a_{11}.
\]
%\textcolor{red}{$\eta_n$ utilise une seule fois!}

The volume grading is defined by the integral linear form
\[
 \ell_n(a)=\sum_{(i,j)\in Q_n}a_{ij},\qquad
 \ell_n(\chi_{D(\lambda)})=|\lambda|.
\]
Thus $\ell_n$ corresponds to the cocharacter $v_n=\sum_{(i,j)\in Q_n}e_{ij}$, which lies in the interior of $\sigma_n$, since $\ell_n$ is positive on $C_{Q_n}\setminus\{0\}$.
With the convention that pullback acts on a function of weight $k$ by multiplication by $\zeta^k$, the volume action is
$\zeta^*x_\lambda=\zeta^{|\lambda|}x_\lambda$ for $\zeta\in\mathbb G_m$.

Write $(A_n)_{s,k}$ for the component of standard degree $s$ and volume $k$, and put $h_n(k)=\dim_\kk(A_n)_k$. The standard basis gives the concrete interpretation
\begin{equation}\label{eq:boundedcount}
 \dim_\kk(A_n)_{s,k}
 =\#\{\pi:|\pi|=k,\ \pi\text{ has }s\text{ nonzero rows},
                  \ \lambda^{(a)}\leq\rho^{(n)}\text{ for every row }a\}.
\end{equation}
In particular $s\leq k\leq d(n)s$, so only finitely many standard degrees contribute to any fixed volume.

For $n\geq2$, set
\[
 X_n=\Proj_{\mathrm{std}}A_n,\qquad \mathcal L_n=\mathcal O_{X_n}(1).
\]
This is a normal projective toric variety of dimension $d(n)-1$, and $\mathcal L_n$ is very ample. It differs from $Y_{L_n}=\Proj R_{L_n}$, whose dimension is $d(n)$: here we projectivize the already dehomogenized algebra $A_n$ using its standard grading.
Normality of the homogeneous ring, and $\dim A_n\geq2$, give weight-compatible isomorphisms
\begin{equation}\label{eq:sections}
 (A_n)_{s,k}\simeq H^0(X_n,\mathcal L_n^{\otimes s})_k,
 \qquad
 h_n(k)=\sum_{s\geq0}\dim_\kk H^0(X_n,\mathcal L_n^{\otimes s})_k.
\end{equation}
The subscript $k$ denotes the weight space for the volume action, with the induced linearization of $\mathcal L_n$.
The relevant integral polytope is
\[
 \rev{\Delta}_n=C_{Q_n}\cap\{a_{11}=1\}
   =\operatorname{conv}\{\chi_{D(\lambda)}:
                    \varnothing\ne\lambda\leq\rho^{(n)}\}.
\]
For $s\geq1$, sections of $\mathcal L_n^{\otimes s}$ have a torus eigenbasis indexed by $s\rev{\Delta}_n\cap M_{Q_n}$; their volume weights are $\ell_n(a)$. At $s=0$ the unique exponent is zero. Thus standard degree records the number of rows, whereas volume records their total weight.

\begin{example}[$U_1$ and the dimension hypothesis]
Here $A_1=\kk[x_{(1)}]$ and both gradings agree. Hence $h_1(k)=1$ for every $k\geq0$, corresponding to a single column of $k$ ones. Although its standard Proj is a point, the dimension hypothesis used to justify \eqref{eq:sections} must not be inferred from normality alone in dimension one; in this example the section identification and weights can simply be checked directly.
\end{example}

\begin{example}[The two gradings on $U_2$]\label{ex:rrquadric}
Recall $A_2=\kk[T,X,Y,Z]/(XY-TZ)$. The volume weights are
$\deg(T,X,Y,Z)=(1,2,2,3)$, while all four standard degrees are one.
Here $X_2\simeq\mathbb P^1\times\mathbb P^1$ and $\mathcal L_2\simeq\mathcal O(1,1)$.
Using $T=u_0v_0$, $X=u_1v_0$, $Y=u_0v_1$, $Z=u_1v_1$, the section
$u_0^{s-a}u_1^av_0^{s-b}v_1^b$ has volume weight $s+a+b$. Therefore
\[
 h_2(k)=\sum_{s\geq0}\#\{(a,b):0\leq a,b\leq s,\ s+a+b=k\}.
\]
For $k=4$, the contributions come from $s=2,3,4$, with dimensions $3,2,1$, respectively. Thus $h_2(4)=6$. These count the plane partitions of volume four whose rows lie in $L_{(2,1)}$; there are thirteen unrestricted plane partitions of that volume.
\end{example}

\subsection{Equivariant Riemann--Roch}
For a finite-dimensional weight-decomposed representation $V=\bigoplus_k V_k$, write
\[\operatorname{ch}_{\mathrm{wt}}(V;q)=\sum_k(\dim_\kk V_k)q^k.\]
This is its character, not its cohomological Chern character. Define
\[
 \chi^{\mathbb G_m}(X_n,\mathcal L_n^{\otimes s};q)
 =\sum_{j\geq0}(-1)^j\operatorname{ch}_{\mathrm{wt}}
       (H^j(X_n,\mathcal L_n^{\otimes s});q).
\]
Demazure vanishing for nef line bundles on complete toric varieties applies even when $X_n$ is singular; see \cite{CLS}. Consequently
\begin{equation}\label{eq:equivrr}
 \Hilb(A_n;q)=\sum_{s\geq0}
     \chi^{\mathbb G_m}(X_n,\mathcal L_n^{\otimes s};q).
\end{equation}
The sum is coefficientwise finite. Equivariant Riemann--Roch computes these Euler characters. In the smooth case, its characteristic-class expression involves the equivariant Chern character of $\mathcal L_n^{\otimes s}$ and the equivariant Todd class of $TX_n$. For singular $X_n$, one uses the singular Riemann--Roch transformation instead of treating $TX_n$ as a vector bundle. Formula~\eqref{eq:equivrr} itself needs neither smoothness nor such a choice of characteristic-class presentation.

For example, on $X_2=\mathbb P^1\times\mathbb P^1$, the character is
$q^s(1+q+\cdots+q^s)^2$. Its value at $q=1$ is $(s+1)^2$, the ordinary dimension of sections. Keeping only this dimension loses the information needed to recover the individual values $h_2(k)$.

\subsection{An integral polytope determined by ordinary partitions}
The polytope $\rev{\Delta}_n$ above is determined by the ordinary partitions in $L_n$: each nonempty partition contributes the zero--one vector of its diagram. In row-by-row coordinates on $\RR^{d(n)}$, it has the explicit inequality description
\[
 \rev{\Delta}_n=\{a:a_{11}=1,\ a_{ij}\geq0,\
 a_{ij}\geq a_{i+1,j},\ a_{ij}\geq a_{i,j+1}\},
\]
where inequalities are included only when their cells belong to $Q_n$. It is an integral polytope of dimension $d(n)-1$; its vertices are exactly the vectors $\chi_{D(\lambda)}$ for nonempty $\lambda\leq\rho^{(n)}$. To see the vertex assertion, a zero--one vector cannot be a nontrivial convex combination of other zero--one vectors.

\begin{theorem}[Lattice points and plane partitions]\label{thm:polytope}
Let $n\geq1$ and $s\geq0$. The lattice points in $s\rev{\Delta}_n$ are in bijection with plane partitions having exactly $s$ nonzero rows, each contained in $\rho^{(n)}$. With $0\rev{\Delta}_n=\{0\}$, one has
\begin{equation}\label{eq:ehrhart}
 \#(s\rev{\Delta}_n\cap\ZZ^{d(n)})
 =\#\{\pi:\pi\text{ has exactly }s\text{ nonzero rows, all in }L_n\}.
\end{equation}
More precisely, for every $k\geq0$,
\begin{align*}
 \#\{a\in s\rev{\Delta}_n\cap\ZZ^{d(n)}:\ell_n(a)=k\}
 &=\#\{\pi:|\pi|=k,\ \pi\text{ has }s\text{ nonzero rows in }L_n\}\\
 &=\dim_\kk(A_n)_{s,k}.
\end{align*}
For $n\geq2$, these numbers are also the dimensions of the weight-$k$ spaces in $H^0(X_n,\mathcal L_n^{\otimes s})$.
\end{theorem}
\begin{proof}
For a plane partition with row chain
$\lambda^{(1)}\geq\cdots\geq\lambda^{(s)}\ne\varnothing$, put
$a=\sum_{h=1}^s\chi_{D(\lambda^{(h)})}$. It is nonnegative, order-reversing, and satisfies $a_{11}=s$, so it lies in $s\rev{\Delta}_n\cap\ZZ^{d(n)}$.
Conversely, given such an $a$, its superlevel ideals
$I_h=\{(i,j):a_{ij}\geq h\}$, $1\leq h\leq s$, are nonempty and nested. They are the diagrams of a unique chain of nonempty partitions in $L_n$. These constructions are inverse, and
$\ell_n(a)=\sum_h|I_h|=|\pi|$.
The statements concerning graded components and sections follow from the standard basis and \eqref{eq:sections}. For $s=0$, both sets consist of the empty object.
\end{proof}

The theorem allows computation in either direction. A lattice-point algorithm applied to the polytope given by the ordinary partitions computes the corresponding plane-partition count. Conversely, a formula for those plane partitions computes the Ehrhart polynomial
$E_n(s)=\#(s\rev{\Delta}_n\cap\ZZ^{d(n)})$. Since $\rev{\Delta}_n$ is integral, $E_n(s)$ is a polynomial of degree $d(n)-1$. Fixing volume as well amounts to imposing the extra linear equation $\ell_n(a)=k$ on the lattice points of $s\rev{\Delta}_n$.

\begin{example}[A square in $\RR^3$]
For $n=2$, the four nonempty ordinary partitions in $L_2$ are
$(1),(2),(1,1),(2,1)$. Their diagram vectors give
\[
 \rev{\Delta}_2=\operatorname{conv}\{(1,0,0),(1,1,0),(1,0,1),(1,1,1)\}
 =\{(1,b,c):0\leq b,c\leq1\}\subseteq\RR^3.
\]
Thus $s\rev{\Delta}_2\cap\ZZ^3=\{(s,b,c):0\leq b,c\leq s\}$ and
$E_2(s)=(s+1)^2$. By Theorem~\ref{thm:polytope}, this is exactly the number of plane partitions with $s$ nonzero rows chosen from $L_{(2,1)}$.
For $s=2$, there are nine. The point $(2,1,1)$ corresponds to the row chain $(2,1)\geq(1)$ and hence to
$\left(\begin{smallmatrix}2&1\\1&0\end{smallmatrix}\right)$.
At volume $k$, the points satisfy $s+b+c=k$. For example, $s=2,k=4$ gives precisely
$(2,0,2),(2,1,1),(2,2,0)$, corresponding respectively to the row chains
$((1,1),(1,1))$, $((2,1),(1))$, and $((2),(2))$.
\end{example}

\begin{example}[A product of triangles in $\RR^5$]
For $n=3$, order the coordinates as
\[(11,12,13,21,31).\]
Then
\[
 \rev{\Delta}_3=\{(1,b,c,d,e):1\geq b\geq c\geq0,\
                                      1\geq d\geq e\geq0\}.
\]
Its nine vertices, corresponding to the hooks $(a,1^b)$ with $1\leq a\leq3$, $0\leq b\leq2$, are
\[
 (1,\mathbf1_{a\geq2},\mathbf1_{a\geq3},
       \mathbf1_{b\geq1},\mathbf1_{b\geq2}).
\]
For each dilation, the two pairs of coordinates can be chosen independently, so
\[
 E_3(s)=\binom{s+2}{2}^{\!2}.
\]
Consequently there are $36$ plane partitions with exactly two nonzero rows in $L_{(3,1,1)}$. This also equals
$h^0(\mathbb P^2\times\mathbb P^2,\mathcal O(2,2))$, in agreement with the standard Proj description of $A_3$.
\end{example}

For volume counts alone, sum the refined equality in Theorem~\ref{thm:polytope} over $s$; only $0\leq s\leq k$ occur. For instance,
\[
 \Hilb(A_2;q)=\frac{1-q^4}{(1-q)(1-q^2)^2(1-q^3)}
 =1+q+3q^2+4q^3+6q^4+8q^5+11q^6+\cdots.
\]
This follows either from the square's weight slices or from the hypersurface equation $XY-TZ$.

More generally, set $B_n=\kk[x_\lambda:\varnothing\ne\lambda\leq\rho^{(n)}]$ with volume grading. A finite graded free resolution of $A_n$ over $B_n$ gives
\[
 \Hilb(A_n;q)=
 \frac{\sum_{j,k}(-1)^j\beta^{(n)}_{j,k}q^k}
      {\prod_{\varnothing\ne\lambda\leq\rho^{(n)}}(1-q^{|\lambda|})},
\]
where $\beta^{(n)}_{j,k}$ are its graded Betti numbers. This provides a computational alternative to lattice-point counting. No complete-intersection assumption is justified by the quadratic form of the generators alone. For instance $U_3$ has nine coordinates and dimension five, but its ideal has nine minimal quadratic minors, so it is not a complete intersection.

\subsection{Brion's formula and counting by volume}\label{subsec:brion}
The lattice-point correspondence also gives access to local counting formulas. Brion's theorem \cite{Brion} expresses the generating polynomial of a lattice polytope as a sum of rational functions associated with its tangent cones at the vertices. We apply it to the polytopes $\rev{\Delta}_n$ of Theorem~\ref{thm:polytope}.

Write $\mathbf z=(z_{ij})_{(i,j)\in Q_n}$ and $\mathbf z^a=\prod_{(i,j)\in Q_n}z_{ij}^{a_{ij}}$. For $s\geq1$, put
\[
 F_{n,s}(\mathbf z)=\sum_{a\in s\rev{\Delta}_n\cap\ZZ^{Q_n}}\mathbf z^a.
\]
For a vertex $v$ of $\rev{\Delta}_n$, let
\[
 T_{n,v}=\cone\{a-v:a\in\rev{\Delta}_n\},\qquad
 \Lambda_n=\{u\in\ZZ^{Q_n}:u_{11}=0\}.
\]
Here $T_{n,v}$ is the cone of tangent directions, rather than the normal cone of the toric fan. Denote by $R_{n,v}(\mathbf z)$ the rational generating function of $T_{n,v}\cap\Lambda_n$. The tangent cone of $s\rev{\Delta}_n$ at $sv$ is $sv+T_{n,v}$: dilation translates its vertex but leaves its cone of directions unchanged.

\begin{proposition}\label{prop:brion}
For $s\geq1$, there is an identity of rational functions
\begin{equation}\label{eq:brion}
 F_{n,s}(\mathbf z)=
 \sum_{v\text{ a vertex of }\rev{\Delta}_n}\mathbf z^{sv}R_{n,v}(\mathbf z).
\end{equation}
Moreover,
\[
 F_{n,s}(q,\ldots,q)
 =\sum_{\substack{\pi\text{ with exactly }s\text{ nonzero rows}\\
                   \text{each row in }L_n}}q^{|\pi|}.
\]
Thus extracting the coefficient of $q^k$ counts these plane partitions of volume $k$, while evaluating at $q=1$ gives $E_n(s)$.
\end{proposition}
\begin{proof}
Apply Brion's theorem in the affine lattice $a_{11}=s$. Every vertex $sv$ is integral, so its cone contributes the monomial $\mathbf z^{sv}$ times the rational generating function of its tangent directions. This proves \eqref{eq:brion}. The second assertion follows from Theorem~\ref{thm:polytope}, since $|\pi|=\sum_{(i,j)\in Q_n}a_{ij}$.
\end{proof}

The specialization in the proposition is always defined on the polynomial $F_{n,s}$. It need not be defined on each cone contribution separately: denominators may vanish when all variables are set equal. One may first combine the rational functions, or substitute $z_{ij}=q\exp(t\eta_{ij})$ with generic integers $\eta_{ij}$ and take the constant term in $t$ after summation. The poles then cancel because the full sum is a polynomial. Likewise, evaluation at $q=1$ is performed after cancellation.

\begin{example}[The four vertex cones of the square]
For $n=2$, use coordinates $(a_{11},a_{12},a_{21})$ and variables $(z_{11},x,y)$. The four vertices of $\rev{\Delta}_2$ are $(1,0,0)$, $(1,1,0)$, $(1,0,1)$, and $(1,1,1)$. Their tangent directions are generated, respectively, by
\[
 \{f_{12},f_{21}\},\quad
 \{-f_{12},f_{21}\},\quad
 \{f_{12},-f_{21}\},\quad
 \{-f_{12},-f_{21}\},
\]
where $f_{12}=(0,1,0)$ and $f_{21}=(0,0,1)$ are the coordinate vectors in $\Lambda_2$. Hence Brion's formula gives
\begin{align*}
 F_{2,s}(z_{11},x,y)=z_{11}^s\bigg(&\frac{1}{(1-x)(1-y)}
 +\frac{x^s}{(1-x^{-1})(1-y)}\\
 &+\frac{y^s}{(1-x)(1-y^{-1})}
 +\frac{x^sy^s}{(1-x^{-1})(1-y^{-1})}\bigg).
\end{align*}
The sum simplifies to
\[
 F_{2,s}(z_{11},x,y)
 =z_{11}^s\frac{1-x^{s+1}}{1-x}\frac{1-y^{s+1}}{1-y}.
\]
Setting $z_{11}=x=y=q$ therefore yields
\[
 F_{2,s}(q,q,q)=q^s(1+q+\cdots+q^s)^2.
\]
For $s=2$, this is $q^2+2q^3+3q^4+2q^5+q^6$: there are nine partitions in total, including three of volume four. This recovers the preceding square example from the four local cone contributions.
\end{example}

Brion's formula applies equally to the integral polytopes constructed recursively in Section~\ref{sec:solid}. It therefore gives finite rational-function expressions for constrained solid- and higher-dimensional partition counts, independently of an unrestricted-volume product formula. The case $s=0$ is treated separately by the empty partition, with generating polynomial $1$.

\subsection{Passage to \texorpdfstring{$U_\infty$}{U-infinity} and MacMahon's series}
The weight forms and standard-degree forms are compatible with extension by zero from $Q_n$ to $Q_{n+1}$. Every plane partition of volume $k$ has all row weights at most $k$. Hence every row lies in $L_n$ as soon as $n\geq k$, and
\begin{equation}\label{eq:stabilize}
 (A_n)_k\simeq(A_\infty)_k\quad(n\geq\max\{1,k\}),\qquad
 h_n(k)=p_2(k),
\end{equation}
where $p_2(k)$ counts unrestricted plane partitions of volume $k$.
Thus the finite toric section counts, equivariant Euler characters, and the refined lattice-point counts of Theorem~\ref{thm:polytope} all stabilize coefficientwise to MacMahon's series:
\[
 \lim_{n\to\infty}\Hilb(A_n;q)
 =\sum_{k\geq0}p_2(k)q^k
 =\prod_{k\geq1}(1-q^k)^{-k}.
\]
For example, $h_2(4)=6$ whereas $h_n(4)=13$ for every $n\geq4$.
The Ehrhart polynomial $E_n(s)$ counts by row number for fixed $n$, not by unrestricted volume. Thus it is distinct from $p_2(k)$. Likewise, we use Riemann--Roch on the finite projective varieties $X_n$, not directly on the infinite-dimensional affine scheme $U_\infty$.

Finally, the main theorem identifies the Hilbert numerator relative to the polynomial algebra on all nonempty ordinary partitions:
\[
 \sum_{\pi\in\PP}S(\pi)q^{|\pi|}
 =\frac{\Hilb(A_\infty;q)}{\Hilb(B;q)}.
\]
This joins the two counting problems: plane partitions index regular functions and graded sections, while their signed incomparability sums encode the numerator arising from the relations among the ordinary-partition generators.

\begingroup
\section{MacMahon's formula from a Grassmannian degeneration}\label{sec:rectangles}
We explain how the classical Grassmannian--Hibi degeneration gives the rectangular plane-partition product, and hence MacMahon's formula. The argument compares the volume Hilbert series of two fibers of a flat family.

\subsection{The rectangular lattice and the Grassmannian}
For $b,c\geq1$, put
\[
L_{b,c}=L_{(c^b)}=\{\lambda:c\geq\lambda_1\geq\cdots\geq\lambda_b\geq0\},
\qquad Q_{b,c}=[b]\times[c].
\]
Let $R_{b,c}=R_{L_{b,c}}$ be the homogeneous Hibi algebra, and set
\[
A_{b,c}=R_{b,c}/(x_\varnothing-1),\qquad \deg x_\lambda=|\lambda|.
\]
Its standard monomial basis consists of chains of nonempty partitions in $L_{b,c}$. Thus its volume Hilbert series counts plane partitions whose nonzero rows have at most $b$ parts, each at most $c$, with no bound on the number of rows. Here $[b]=\{1,\ldots,b\}$ and $A_{b,c}=A_{Q_{b,c}}$ in the notation of Section~2.

The varieties $U_{b,c}=\Spec A_{b,c}$ have dimension $bc$ and give the same affine inverse limit $U_\infty$ as the weight-generated approximations: the inclusions $L_n\subseteq L_{n,n}$ and $L_{b,c}\subseteq L_{bc}$ make the two families cofinal. The notation $L_{bc}$ here means the weight-generated lattice with cutoff $bc$, whereas $L_{b,c}$ denotes the rectangular lattice. For example, $Q_4$ has eight cells, while the square $Q_{4,4}$ has sixteen. Rectangles give larger models for a fixed volume cutoff, but their Hilbert series admit the product evaluation below.

The Grassmannian $\operatorname{Gr}(b,b+c)$ parametrizes $b$-dimensional subspaces of $\kk^{b+c}$. Its Pl\"ucker coordinate $p_I$, for $I=\{i_1<\cdots<i_b\}$, is the corresponding maximal minor of a matrix representing the subspace. The indexing
\[
I(\lambda)=\{r+\lambda_{b+1-r}:1\leq r\leq b\}
\]
identifies the lattice of these subsets with $L_{b,c}$. Write $p_\lambda=p_{I(\lambda)}$, and let
\[
S_{b,c}=\kk[p_\lambda:\lambda\in L_{b,c}]/I_{\mathrm{Pl}}
\]
be the homogeneous coordinate ring, where $I_{\mathrm{Pl}}$ is the Pl\"ucker ideal. In addition to standard degree $\deg_{\mathrm{std}}p_\lambda=1$, it has the volume grading
\begin{equation}\label{eq:grass-volume}
\deg_{\mathrm{vol}}p_I=\sum_{r=1}^b(i_r-r)=|\lambda|.
\end{equation}
Indeed, Pl\"ucker relations preserve the multiset of column indices and standard degree, so they are homogeneous for this grading. In particular, $p_\varnothing=p_{\{1,\ldots,b\}}$ has volume degree zero.

\subsection{The deformation and its flatness}
We use maximum-weight initial forms and fix the following integral weight on the Pl\"ucker variables:
\begin{equation}\label{eq:explicit-grass-weight}
\omega_\lambda=\sum_{i=1}^b(b-i)\lambda_i,
\qquad\text{equivalently}\qquad
\omega(p_I)=\sum_{r=1}^b(r-1)(i_r-r),
\end{equation}
where $I=\{i_1<\cdots<i_b\}=I(\lambda)$. In particular, $\omega_\varnothing=0$. The weight of $p^\alpha=\prod_\lambda p_\lambda^{\alpha_\lambda}$ is
$\omega\cdot\alpha=\sum_\lambda\omega_\lambda\alpha_\lambda$.
With this choice, the diagonal Grassmannian--Hibi degeneration \cite{GL96,GL97,FM} gives
\[
\operatorname{in}_\omega(I_{\mathrm{Pl}})=J_{b,c},
\]
where $J_{b,c}$ is generated by
$p_\lambda p_\mu-p_{\lambda\wedge\mu}p_{\lambda\vee\mu}$ for incomparable indices.

To see the origin of the weight, represent the Pl\"ucker coordinates by the maximal minors of a generic matrix $(a_{rs})$ and give $a_{rs}$ weight $(r-1)s$. The strict rearrangement inequality makes $\prod_r a_{r,i_r}$ the unique maximum-weight term of the minor $p_I$. Its weight is $\sum_r(r-1)i_r$; subtracting the constant $\sum_r(r-1)r$ gives \eqref{eq:explicit-grass-weight}. Products indexed by chains have distinct diagonal monomials: sorting the column indices separately in each matrix row recovers the unique chain. The diagonal monomial algebra therefore has exactly the Hibi relations. The standard monomial basis of the Grassmannian then identifies it with the associated graded algebra, proving the stated initial-ideal equality.

\paragraph{Remark on the choice of weight.}
The weight is not unique. More generally, strictly increasing integers $\theta_1<\cdots<\theta_b$ give the same diagonal degeneration with
\[
\omega^{(\theta)}(p_I)=\sum_{r=1}^b\theta_r(i_r-r).
\]
The preceding diagonal argument applies unchanged. Our choice is $\theta_r=r-1$. One can also replace $\omega_\lambda$ by $m\omega_\lambda+\ell|\lambda|$, with $m$ a positive integer and $\ell$ any integer: the volume term is constant on each volume-homogeneous relation, so the initial ideal remains $J_{b,c}$. The fibers and volume Hilbert series are unchanged, although the powers of the deformation parameter can change; multiplication by $m$ corresponds to the base change $t\mapsto t^m$. More generally, any integral weight with the same initial ideal works, after normalizing its empty-partition weight to zero. An arbitrary weight need not have this property.

For $b=2$, formula~\eqref{eq:explicit-grass-weight} reads $\omega(p_{ij})=j-2$, or $\omega_\lambda=\lambda_1$. Thus, for $i<j<k<l$, the relation
\[
p_{il}p_{jk}-p_{ik}p_{jl}+p_{ij}p_{kl}=0
\]
has its first two terms of equal maximum weight, exceeding the weight of the third by $k-j$. Its homogenization is
\[
p_{il}p_{jk}-p_{ik}p_{jl}+t^{k-j}p_{ij}p_{kl}=0.
\]
In particular, the case $(i,j,k,l)=(1,2,3,4)$ gives the family $XY-TZ+tW$ used below.

Introduce a parameter $t$. For a nonzero $f=\sum_\alpha c_\alpha p^\alpha\in I_{\mathrm{Pl}}$, put
\[
d_\omega(f)=\max_{c_\alpha\ne0}\omega\cdot\alpha,
\qquad
\widetilde f=\sum_\alpha c_\alpha t^{d_\omega(f)-\omega\cdot\alpha}p^\alpha.
\]
Let $\widetilde I$ be generated by these homogenizations for all $f\in I_{\mathrm{Pl}}$, and define
\[
\mathcal R_{b,c}=\kk[t,p_\lambda:\lambda\in L_{b,c}]/\widetilde I,
\qquad
\mathcal A_{b,c}=\mathcal R_{b,c}/(p_\varnothing-1).
\]
Homogenizing all relations, or a suitable weight Gr\"obner basis, is essential: arbitrary perturbations of the Hibi relations need not give a flat family. The weight $\omega$ controls the degeneration; it is distinct from the volume grading. Both standard degree and volume extend to $\mathcal R_{b,c}$ by assigning $t$ degree zero.

\begin{proposition}\label{prop:grass-family}
The algebra $\mathcal A_{b,c}$ is a flat volume-graded $\kk[t]$-algebra. Its fiber at $t=0$ is $A_{b,c}$, and its fiber at any $t=a\ne0$ is the coordinate ring of the big cell of $\operatorname{Gr}(b,b+c)$. Each volume component is a finite free $\kk[t]$-module.
\end{proposition}
\begin{proof}
The weight-homogenization construction realizes $\mathcal R_{b,c}$ as the subalgebra
\[
\kk[t,t^{\omega_\lambda}\overline p_\lambda:\lambda\in L_{b,c}]
\subseteq S_{b,c}[t,t^{-1}],
\]
where $\overline p_\lambda$ is the class of $p_\lambda$ in $S_{b,c}$. The kernel of the displayed substitution is exactly $\widetilde I$. Since $S_{b,c}$ is a domain, this realizes $\mathcal R_{b,c}$ inside a domain and embeds $\kk[t]$ in it. It is therefore torsion-free over $\kk[t]$.

To justify the affine-chart step, use standard degree:
\[
\mathcal A_{b,c}\simeq
\bigl(\mathcal R_{b,c}[p_\varnothing^{-1}]\bigr)_{\mathrm{std},0}.
\]
This is a subalgebra of a localization of the same domain, so it too is torsion-free over $\kk[t]$. A module over the principal ideal domain $\kk[t]$ is flat if and only if it is torsion-free; thus $\mathcal A_{b,c}$ is flat.

At $t=0$, the relations specialize to the Hibi ideal. For $t=a\ne0$, the substitution $p_\lambda\mapsto a^{\omega_\lambda}\overline p_\lambda$ identifies the relations with the Pl\"ucker relations; it preserves $p_\varnothing$ because $\omega_\varnothing=0$. The fibers are consequently $A_{b,c}$ and $S_{b,c}/(p_\varnothing-1)$, respectively.

Finally, after imposing $p_\varnothing=1$, all remaining generators have strictly positive volume degrees. Every volume component is therefore finite over $\kk[t]$. It is torsion-free, hence free over this principal ideal domain. Its rank is the dimension of the corresponding component in every fiber.
\end{proof}

\subsection{The weights on the polynomial fiber}
The big cell $p_\varnothing\ne0$ consists of subspaces with a unique matrix representative
\[
[\,I_b\mid Y\,],\qquad Y=(y_{rj})_{1\leq r\leq b,\ 1\leq j\leq c}.
\]
Its coordinate ring is the polynomial ring $\kk[y_{rj}]$. Up to sign, $y_{rj}$ is the Pl\"ucker ratio obtained by replacing column $r$ of $\{1,\ldots,b\}$ by column $b+j$. Formula~\eqref{eq:grass-volume} gives
\[
\deg y_{rj}=b+j-r.
\]
With $x_{ij}=y_{b+1-i,j}$, the graded polynomial fiber is therefore
\begin{equation}\label{eq:grass-polynomial}
B_{b,c}=\kk[x_{ij}:1\leq i\leq b,\ 1\leq j\leq c],
\qquad \deg x_{ij}=i+j-1.
\end{equation}
These are the degrees inherited from volume, not the auxiliary degeneration weights.

\begin{example}[The rectangle $2\times2$]
For the partitions $(1),(2),(1,1),(2,1),(2,2)$, write the corresponding generators as $T,X,Y,Z,W$, of volume degrees $1,2,2,3,4$. Then
\[
A_{2,2}=\kk[T,X,Y,Z,W]/(XY-TZ).
\]
The affine family is
\[
\mathcal A=\kk[t,T,X,Y,Z,W]/(XY-TZ+tW),\qquad \deg t=0.
\]
It comes from the Pl\"ucker relation for $\operatorname{Gr}(2,4)$ after setting $p_{12}=1$. At $t=0$, it gives $A_{2,2}$. At $t=a\ne0$, eliminating $W=(TZ-XY)/a$ gives $\kk[T,X,Y,Z]$, with weights $1,2,2,3$. The defining polynomial is irreducible, since as a polynomial in $W$ its coefficients $t$ and $XY-TZ$ are coprime. Thus the quotient is a domain, contains $\kk[t]$, and is torsion-free over $\kk[t]$; this also verifies flatness directly. The resulting Hilbert series is
\[
\Hilb(A_{2,2};q)=\frac{1}{(1-q)(1-q^2)^2(1-q^3)}.
\]
\end{example}

\subsection{MacMahon's formula and the graph numerator}
\begin{theorem}\label{thm:grass-macmahon}
Let $\mathcal P_{b,c}$ be the set of plane partitions whose nonzero rows belong to $L_{b,c}$. Then
\[
\sum_{\pi\in\mathcal P_{b,c}}q^{|\pi|}
=\Hilb(A_{b,c};q)
=\prod_{i=1}^b\prod_{j=1}^c(1-q^{i+j-1})^{-1}.
\]
Consequently,
\[
\sum_{\pi\in\PP}q^{|\pi|}=\prod_{m\geq1}(1-q^m)^{-m}.
\]
\end{theorem}
\begin{proof}
The chain basis gives the first equality. Proposition~\ref{prop:grass-family} identifies the dimensions of each volume component of $A_{b,c}$ and $B_{b,c}$. Summing the monomials of the polynomial ring \eqref{eq:grass-polynomial} gives the finite product.

For the limit, fix $N$. Every plane partition of volume at most $N$ lies in $\mathcal P_{b,c}$ once $b,c\geq N$. For $m\leq N$, all $m$ pairs of positive integers with $i+j-1=m$ also lie in that rectangle. Both sides therefore stabilize coefficientwise to the unrestricted identity. No compatibility between the different finite deformation families is needed for this argument.
\end{proof}

\begin{corollary}[Rectangular version of Theorem~\ref{thm:main}]\label{cor:rectangular-graph}
Let $b,c\geq1$, let $\mathcal P_{b,c}$ be the plane partitions whose nonzero rows belong to $L_{b,c}$, with no bound on the number of rows, and retain the statistic $S$ of Definition~\ref{def:defect}, including $S(\varnothing)=1$. Let $G_{b,c}$ be the incomparability graph on $V_{b,c}=L_{b,c}\setminus\{\varnothing\}$ and set $u_\lambda=q^{|\lambda|}$. Then
\begin{equation}\label{eq:rectangular-defect}
\sum_{\pi\in\mathcal P_{b,c}}S(\pi)q^{|\pi|}
=\mathcal S_{G_{b,c}}(\mathbf u,-1)
=\frac{\displaystyle\prod_{\lambda\in V_{b,c}}(1-q^{|\lambda|})}
{\displaystyle\prod_{i=1}^b\prod_{j=1}^c(1-q^{i+j-1})}.
\end{equation}
\end{corollary}
\begin{proof}
The edge-ideal quotient has the same chain basis as $A_{b,c}$, hence the same volume Hilbert series. The graph expansion of Section~\ref{sec:proof} therefore gives
\[
\mathcal S_{G_{b,c}}(\mathbf u,-1)
=\Hilb(A_{b,c};q)\prod_{\lambda\in V_{b,c}}(1-q^{|\lambda|}).
\]
Theorem~\ref{thm:grass-macmahon} evaluates the Hilbert series and proves the second equality in \eqref{eq:rectangular-defect}.

For the first equality, group the finite subgraphs without isolated vertices by the plane partition obtained by sorting the columns of their vertex partitions. Sorting preserves the bounds of at most $b$ entries per row and entries at most $c$. Conversely, every matrix in $\mathcal M(\pi)$ for $\pi\in\mathcal P_{b,c}$ has distinct nonzero rows in $L_{b,c}$: column permutations preserve these bounds, and its first column has only positive entries. Its selections of incomparable row pairs covering every row are exactly the edge sets of subgraphs with this vertex set and no isolated vertices. As in the proof of Theorem~\ref{thm:main}, their signs are $(-1)$ to the number of selected pairs, so their sum is $S(\pi)$. Sorting also preserves volume. The empty graph supplies $S(\varnothing)=1$.
\end{proof}

\subsection{Toric gradings and the refined signed numerator}\label{subsec:toric-refinement}
The monomial realization of Proposition~\ref{prop:cone} identifies
$A_{b,c}$ with $\kk[C_{b,c}\cap\ZZ^{bc}]$, where $C_{b,c}=C_{Q_{b,c}}=\sigma_{Q_{b,c}}^\vee$ is the cone of nonnegative arrays decreasing in both grid directions. Under this identification, a standard monomial with rows $\lambda^{(1)},\ldots,\lambda^{(r)}$ has exponent
\[
h=\sum_{a=1}^r\chi_{D(\lambda^{(a)})},\qquad
h_{11}=r,\qquad \ell(h):=\sum_{i,j}h_{ij}=|\pi|.
\]
Thus the toric interpretation already gives the bijection between plane partitions and integer points of the cone. Retaining both gradings, put
\[
F_{b,c}(t,q)=\sum_{h\in C_{b,c}\cap\ZZ^{bc}}t^{h_{11}}q^{\ell(h)},
\qquad
N_{b,c}(t,q)=\sum_{\pi\in\mathcal P_{b,c}}S(\pi)t^{r(\pi)}q^{|\pi|},
\]
where $r(\pi)$ is the number of nonzero rows. Here $t$ records standard degree, independently of the deformation parameter used above. In particular, $F_{b,c}$ is the bigraded Hilbert series of $A_{b,c}$.

Label the cells by $\nu(i,j)=(i-1)c+j$. A linear extension is a listing of these labels respecting the grid order; write $\operatorname{LinExt}(Q_{b,c})$ for their set. For such a word $w=w_1\cdots w_{bc}$, define the descent set, descent number, and major index by
\[
\operatorname{Des}(w)=\{k:w_k>w_{k+1}\},\quad
\operatorname{des}(w)=|\operatorname{Des}(w)|,\quad
\operatorname{maj}(w)=\sum_{k\in\operatorname{Des}(w)}k.
\]
These are the usual permutation statistics \cite[\S1.4]{StanleyEC1}. We denote their enumerator by
\[
E_{b,c}(t,q)=\sum_{w\in\operatorname{LinExt}(Q_{b,c})}
 t^{\operatorname{des}(w)}q^{\operatorname{maj}(w)}.
\]
The letter $E$ distinguishes this Eulerian-type polynomial from the algebra $A_{b,c}$.

\begin{proposition}\label{prop:refined-toric}
Set $p_{b,c}(k)=\#\{\lambda\in L_{b,c}:|\lambda|=k\}$. Then
\begin{align}
F_{b,c}(t,q)
&=\frac{N_{b,c}(t,q)}{\prod_{\lambda\in V_{b,c}}(1-tq^{|\lambda|})}
=\frac{E_{b,c}(t,q)}{\prod_{k=1}^{bc}(1-tq^k)},\label{eq:refined-toric}\\
N_{b,c}(t,q)
&=E_{b,c}(t,q)\prod_{k=1}^{bc}(1-tq^k)^{p_{b,c}(k)-1}.\label{eq:refined-factor}
\end{align}
Thus the signed numerator factors through a polynomial with nonnegative coefficients. The coefficient $[t^s q^m]F_{b,c}$ counts the integer points with $h_{11}=s$ and $\ell(h)=m$, equivalently plane partitions of volume $m$ with exactly $s$ rows in $L_{b,c}$.
\end{proposition}
\begin{proof}
The first equality in \eqref{eq:refined-toric} follows from the graph expansion with $u_\lambda=tq^{|\lambda|}$: the grouping in Corollary~\ref{cor:rectangular-graph} preserves both the number of rows and volume.

For the second, use the classical linear-extension decomposition of $P$-partitions \cite[\S3.15, Lemma~3.15.3 and Theorem~3.15.5]{StanleyEC1}, whose geometric counterpart is the canonical triangulation of order polytopes \cite[\S5]{StanleyPolytopes}. Order the coordinates of $h$ decreasingly, breaking ties by increasing label. This assigns $h$ to a unique extension $w$ and partitions the exponent cone into half-open unimodular cones. On the cone indexed by $w$, set
$d_k=h_{w_k}-h_{w_{k+1}}$, with $h_{w_{bc+1}}=0$. Then $d_k\geq0$, with $d_k\geq1$ required at descents. Since $(1,1)$ is first in every extension,
\[
h_{11}=\sum_k d_k,\qquad \ell(h)=\sum_k k d_k.
\]
Its lattice-point series is therefore
$t^{\operatorname{des}(w)}q^{\operatorname{maj}(w)}/\prod_k(1-tq^k)$.
Summing gives \eqref{eq:refined-toric}, and comparison of the denominators gives \eqref{eq:refined-factor}. Every weight $1\leq k\leq bc$ occurs in $L_{b,c}$, so all exponents $p_{b,c}(k)-1$ are nonnegative.
\end{proof}

\begin{example}
For the $2\times2$ grid, labeled row by row, the extensions are $1234$ and $1324$. Their descent sets are $\varnothing$ and $\{2\}$, so
\[
E_{2,2}(t,q)=1+tq^2,\qquad N_{2,2}(t,q)=1-t^2q^4,
\]
\[
F_{2,2}(t,q)=\frac{1+tq^2}{(1-tq)(1-tq^2)(1-tq^3)(1-tq^4)}.
\]
The negative term corresponds to the single edge joining $(2)$ and $(1,1)$, whose column sorting gives
$\pi=\left(\begin{smallmatrix}2&1\\1&0\end{smallmatrix}\right)$ with $S(\pi)=-1$. Specializing $t=1$ recovers the rectangular MacMahon product and its signed numerator.
\end{example}

\begin{corollary}\label{cor:refined-limit}
Let $F(t,q)=\sum_{\pi\in\PP}t^{r(\pi)}q^{|\pi|}$ and
$N(t,q)=\sum_{\pi\in\PP}S(\pi)t^{r(\pi)}q^{|\pi|}$.
The coefficientwise limit $E_\infty(t,q)=\lim_{b,c\to\infty}E_{b,c}(t,q)$ exists in $\ZZ_{\geq0}[t][[q]]$, and
\[
F(t,q)=\frac{E_\infty(t,q)}{\prod_{k\geq1}(1-tq^k)},\qquad
N(t,q)=E_\infty(t,q)\prod_{k\geq1}(1-tq^k)^{p(k)-1}.
\]
At $t=1$, one has $E_\infty(1,q)=\prod_{k\geq1}(1-q^k)^{-(k-1)}$, and the second identity becomes Theorem~\ref{thm:main}.
\end{corollary}
\begin{proof}
For a fixed volume bound $M$, all plane partitions and ordinary partitions involved occur once $b,c\geq M$. The products in \eqref{eq:refined-toric} stabilize through $q^M$, hence so does $E_{b,c}$. Its coefficients are nonnegative. Passing to the limit proves the identities, and Theorem~\ref{thm:grass-macmahon} evaluates the specialization at $t=1$.
\end{proof}
The factorization links the signed incomparability expansion to the positive numerator supplied by decomposing the toric exponent cone. Proposition~\ref{prop:higher-refined} extends this factorization to the finite toric approximations in every dimension.

\par\endgroup

\section{From the lattice of plane partitions to solid partitions}\label{sec:solid}
The construction extends by one dimension. The distributive lattice is now the set of plane partitions, and chains in this lattice are solid partitions. The toric, Gr\"obner, and polytope constructions remain valid; the final evaluation of the volume generating series is the step that no longer follows from MacMahon's formula.

\editnov{There are two parallel choices of finite approximations, extending $L_n$ and $L_{b,c}$, respectively. The first is the lattice generated by lower-dimensional partitions of volume at most $n$; its elements need not themselves have volume at most $n$. The second consists of all lower-dimensional partitions whose diagrams lie in a fixed product box $[b_1]\times\cdots\times[b_r]$, a hypercube when all side lengths agree. Both choices use lower-dimensional partitions as generators for the next-dimensional toric model. We develop the volume-generated family first, while the box-bounded family provides the rectangular counterpart and is used in Section~\ref{subsec:higher-refinement}. The two families are cofinal and define the same direct-limit algebra and inverse-limit affine scheme, as explained in Section~\ref{subsec:all-dimensions}.}

\subsection{The lattice of plane partitions and its finite approximations}
Order $\PP$ componentwise, and include the empty array. Meet and join are entrywise minimum and maximum, so $\PP$ is a distributive lattice and
\[
 |\pi|+|\tau|=|\pi\wedge\tau|+|\pi\vee\tau|.
\]
Its diagram realization is the set of finite order ideals in $\ZZ_{>0}^3$:
\[
 D_3(\pi)=\{(i,j,k):1\leq k\leq\pi_{ij}\}.
\]
The nonzero join-irreducibles are precisely the rectangular plane partitions
$R_{a,b,c}$, equal to $c$ on $\{1,\ldots,a\}\times\{1,\ldots,b\}$ and zero elsewhere. Their diagrams are the principal boxes of the product poset, and their volumes are $abc$.
Every plane partition is the unique irredundant join of the rectangles associated with the maximal cells of its diagram, by the same ideal argument as Proposition~\ref{prop:rect}.

Let $L_n^{(3)}$ be the sublattice generated by plane partitions of volume at most $n$, and set
\[
 \rho^{(n,3)}_{ij}=\left\lfloor\frac{n}{ij}\right\rfloor,
 \qquad Q_n^{(3)}=\{(i,j,k)\in\ZZ_{>0}^3:ijk\leq n\}.
\]
Here the superscript $3$ refers to the dimension of the diagram poset.
\begin{proposition}\label{prop:solidfinite}
One has $L_n^{(3)}=\{\pi:\pi\leq\rho^{(n,3)}\}\simeq\Id(Q_n^{(3)})$.
In particular, with
\[
 d_3(n)=|Q_n^{(3)}|
 =\sum_{i,j\geq1}\left\lfloor\frac{n}{ij}\right\rfloor
 =\sum_{m=1}^n\tau_3(m)\sim \frac12 n(\log  n)^2,
\]
where $\tau_3(m)=\#\{(a,b,c)\in\ZZ_{>0}^3:abc=m\}$, the normal affine toric variety
$U_n^{(3)}=\Spec A_n^{(3)}$, $A_n^{(3)}:=A_{Q_n^{(3)}}$, has dimension $d_3(n)$.
\end{proposition}
\begin{proof}
A plane partition of volume at most $n$ satisfies $ij\pi_{ij}\leq n$, hence lies below $\rho^{(n,3)}$. Conversely, every principal box contained in $Q_n^{(3)}$ has volume at most $n$ and is among the generators. Their joins give every ideal of $Q_n^{(3)}$. The dimension follows from the dehomogenized Hibi construction.

We approximate the number of positive integer lattice points under the hyperbolic surface \(xyz=n\) by the volume of the corresponding region in \(\mathbb{R}_{\ge1}^3\), or equivalently, sum over \(k\) first:
\[
d_3(n)=\sum_{ij\le n}\Bigl\lfloor\frac{n}{ij}\Bigr\rfloor
\sim\iint_{\substack{xy\le n\\x,y\ge1}}\frac{n}{xy}\,dx\,dy.
\]
Evaluating this double integral over the region \(1\le x\le n\) and \(1\le y\le\frac{n}{x}\) gives
\begin{align*}
\int_1^n\Biggl(\int_1^{n/x}\frac{n}{xy}\,dy\Biggr)dx
&=\int_1^n\frac{n}{x}\log \Bigl(\frac{n}{x}\Bigr)\,dx
=n\int_1^n\frac{\log  n-\log  x}{x}\,dx\\
&=n\Biggl[(\log  n)\log  x-\frac12(\log  x)^2\Biggr]_1^n
\\ &=n\Biggl((\log  n)^2-\frac12(\log  n)^2\Biggr)
=\frac12 n(\log  n)^2.
\end{align*}

\end{proof}
The sum over $i,j$ has only finitely many nonzero terms. In coordinates indexed by $Q_n^{(3)}$, the exponent cone and its dual are
\begin{align*}
 C_n^{(3)}&=\{a\in\RR_{\geq0}^{Q_n^{(3)}}:a_u\geq a_v\text{ whenever }u\leq v\},\\
 \sigma_n^{(3)}&=\cone\{e_u-e_v:u\lessdot v;\ e_u:u\text{ maximal}\}.
\end{align*}
Covers increase one grid coordinate by one. A cell $(i,j,k)$ is maximal exactly when
$(i+1)jk$, $i(j+1)k$, and $ij(k+1)$ all exceed $n$.

\begin{example}[The first branching example]
For $n=1$, $U_1^{(3)}\simeq\mathbb A^1$.
For $n=2$, the cells, in order, are $111,211,121,112$, where $ijk$ abbreviates the triple $(i,j,k)$. The first is below the other three, which are incomparable. Thus $d_3(2)=4$ and
\[
 C_2^{(3)}=\{(a,b,c,d):a\geq b,c,d\geq0\},
\]
\[
 \sigma_2^{(3)}=\cone\{e_2,e_3,e_4,e_1-e_2,e_1-e_3,e_1-e_4\}\subseteq\RR^4.
\]
There are eight nonempty ideals. Their indicator vectors are $(1,\varepsilon_1,\varepsilon_2,\varepsilon_3)$ with $\varepsilon_i\in\{0,1\}$. Hence $U_2^{(3)}$ is the affine cone over the Segre embedding of $(\mathbb P^1)^3$, and its standard Proj has dimension three.
\end{example}

\editnov{For comparison, the box-bounded lattice of plane partitions is
\[
 L_{(a,b,c)}^{(3)}=\{\pi:D_3(\pi)\subseteq[a]\times[b]\times[c]\}
 =\{\pi:\pi\leq R_{a,b,c}\}\simeq\Id([a]\times[b]\times[c]).
\]
Its affine Hibi variety has dimension $abc$, and its standard monomials correspond to solid partitions with an arbitrary number of slices, each bounded by $R_{a,b,c}$. For instance, the box $[2]^3$ has eight cells and gives an eight-dimensional affine variety; it differs from the four-cell poset $Q_2^{(3)}$ in the preceding example. The hypercube polytope appearing there is a slice of the exponent cone, not a hypercube of diagram cells.}

\subsection{Solid partitions as functions and lattice points}
A solid partition is a finitely supported array $\Pi=(\Pi_{ijk})_{i,j,k\geq1}$ of nonnegative integers, weakly decreasing in each coordinate. Write $\mathcal{SP}$ for their set, $|\Pi|=\sum_{i,j,k}\Pi_{ijk}$, and $p_3(m)=\#\{\Pi:|\Pi|=m\}$.
Its nonzero slices $\pi^{(h)}=(\Pi_{hjk})_{j,k\geq1}$ form a decreasing chain of nonempty plane partitions. Define
\[
 A_\infty^{(3)}=\varinjlim_n A_n^{(3)},\qquad
 U_\infty^{(3)}=\Spec A_\infty^{(3)}=\varprojlim_n U_n^{(3)}.
\]
As before, the algebra inclusions are induced by extension of ideal-indicator vectors by zero. Equivalently,
\[
 A_\infty^{(3)}=
 \kk[y_\pi:\pi\in\PP]\big/
 (y_\varnothing-1,\ y_\pi y_\tau-y_{\pi\wedge\tau}y_{\pi\vee\tau}:
                       \pi\parallel\tau).
\]
The standard monomials $y_{\pi^{(1)}}\cdots y_{\pi^{(s)}}$ identify solid partitions with a basis of the global functions on $U_\infty^{(3)}$. Giving each nonempty generator $y_\pi$ volume degree $|\pi|$, this proves
\begin{equation}\label{eq:solidseries}
 \Hilb(A_\infty^{(3)};q)=F_3(q):=\sum_{m\geq0}p_3(m)q^m.
\end{equation}
This equality does not require an evaluation of $F_3$.

The standard-degree and volume forms on $C_n^{(3)}$ are respectively
$a\mapsto a_{111}$ and $\ell_n^{(3)}(a)=\sum_{u\in Q_n^{(3)}}a_u$. Put
\[
 \rev{\Delta}_n^{(3)}=C_n^{(3)}\cap\{a_{111}=1\}
 =\operatorname{conv}\{\chi_{D_3(\pi)}:
                   \varnothing\ne\pi\leq\rho^{(n,3)}\}.
\]
\begin{theorem}\label{thm:solidpolytope}
For $s\geq0$, the lattice points of $s\rev{\Delta}_n^{(3)}$ correspond bijectively to solid partitions with exactly $s$ nonzero slices, each in $L_n^{(3)}$. Under this correspondence, volume equals $\ell_n^{(3)}(a)$.
Consequently the unrestricted-volume count with $s$ slices is the Ehrhart polynomial of the integral polytope $\rev{\Delta}_n^{(3)}$, of degree $d_3(n)-1$.
\end{theorem}
\begin{proof}
Map a slice chain to the sum of its diagram indicators. Recover the chain by the superlevel ideals of the resulting order-reversing integer vector. This is the proof of Theorem~\ref{thm:polytope} on the three-dimensional product poset.
\end{proof}
For $n\geq2$, the same count is the dimension of
$H^0(\Proj_{\mathrm{std}}A_n^{(3)},\mathcal O(s))$, and fixing volume gives a weight space. Normality, toric vanishing, and the equivariant Riemann--Roch interpretation of Section~\ref{sec:geometry} apply without change.
For fixed volume $m$, all slices have volume at most $m$, so the degree-$m$ components stabilize to $(A_\infty^{(3)})_m$ once $n\geq\max\{1,m\}$.

\begin{example}[A cube and an exact solid-partition count]
For $n=2$,
$\rev{\Delta}_2^{(3)}=\{(1,b,c,d):0\leq b,c,d\leq1\}\subseteq\RR^4$.
Therefore there are exactly $(s+1)^3$ solid partitions with $s$ nonzero slices in $L_2^{(3)}$. For $s=2$, there are $27$.
For example, the point $(2,1,1,0)$ has superlevel ideals
$\{111,211,121\}$ and $\{111\}$; these give the slices
$\left(\begin{smallmatrix}1&1\\1&0\end{smallmatrix}\right)$ and $(1)$, of total volume four. More generally the weight character of degree $s$ is
$q^s(1+q+\cdots+q^s)^3$.
\end{example}

\subsection{The signed identity and the missing product evaluation}
Let $G_\infty^{(3)}$ be the incomparability graph on nonempty plane partitions, and let
$B^{(3)}=\kk[y_\pi:\varnothing\ne\pi\in\PP]$, with volume grading. There are $p_2(m)$ generators of weight $m$, hence
\[
 \Hilb(B^{(3)};q)=\prod_{m\geq1}(1-q^m)^{-p_2(m)}.
\]
Hibi's Gr\"obner basis theorem and the subgraph formula of \cite{MM} give
\begin{equation}\label{eq:solidgraph}
 \sum_{H\in\mathcal H(G_\infty^{(3)})}
 (-1)^{|E(H)|}q^{\sum_{\pi\in V(H)}|\pi|}
 =F_3(q)\prod_{m\geq1}(1-q^m)^{p_2(m)}.
\end{equation}
All limits remain coefficientwise finite.

This numerator can also be grouped by solid partitions. For $\Pi\ne\varnothing$, let $\mathcal M_3(\Pi)$ consist of arrays obtained by permuting entries along the first coordinate independently for each fixed pair $(j,k)$, requiring that their nonzero slices be distinct plane partitions arranged in decreasing lexicographic order. For this order, read each slice on the common finite rectangular support row by row. Let $\rev{\mathcal E}_3(A)$ be the incomparable pairs of slices, under componentwise comparison, and define
\[
 S_3(\Pi)=\sum_{A\in\mathcal M_3(\Pi)}
 \sum_{\substack{D\subseteq\rev{\mathcal E}_3(A)\\D\text{ covers every slice}}}(-1)^{|D|},
 \qquad S_3(\varnothing)=1.
\]
Sorting the first-coordinate fibers of a set of plane partitions produces a solid partition: order statistics preserve the inequalities in both remaining coordinates. The canonical slice order makes each vertex set occur once. Thus the same fiber argument as in Section~\ref{sec:proof} proves
\begin{equation}\label{eq:solididentity}
 {\displaystyle
 \sum_{\Pi\in\mathcal{SP}}S_3(\Pi)q^{|\Pi|}
 =F_3(q)\prod_{m\geq1}(1-q^m)^{p_2(m)}.}
\end{equation}
This is an exact identity, but not an explicit product evaluation of its left-hand side: the unknown solid-partition volume series $F_3$ remains.

The step used for plane partitions was to substitute MacMahon's formula for $\Hilb(A_\infty;q)$. The analogous proposed product
\[
 M_3(q)=\prod_{m\geq1}(1-q^m)^{-\binom{m+1}{2}}
\]
cannot be substituted for $F_3(q)$. Already
\[
 [q^6]F_3(q)=140,\qquad [q^6]M_3(q)=141;
\]
see \cite{Atkin}. Accordingly, replacing the right-hand side of \eqref{eq:solididentity} by
$\prod_{m\geq1}(1-q^m)^{p_2(m)-\binom{m+1}{2}}$ would be false.
Amanov and Yeliussizov \cite{AY} recover $M_3$ with corner-hook volume, a different statistic; that result does not evaluate the ordinary-volume Hilbert series used here.

The obstruction is therefore confined to the final explicit evaluation. The normal toric approximations, their dimensions and cones, the solid-partition basis, the finite Ehrhart and section counts, and the graph expansion are all available. They provide exact finite computations and identities relating $F_3$ to a signed numerator, but do not by themselves supply a MacMahon-type formula for ordinary solid-partition volume.

\subsection{The same construction in every dimension}\label{subsec:all-dimensions}
There is no structural reason to stop at solid partitions. Fix $d\geq1$, and let $\mathcal P^{(d)}$ be the lattice of finitely supported nonnegative arrays on $\ZZ_{>0}^d$ that decrease in each coordinate. Write $p_d(m)$ for their volume counts and
$F_d(q)=\sum_{m\geq0}p_d(m)q^m$. Thus $\mathcal P^{(1)}=\Par$, $\mathcal P^{(2)}=\PP$, and $\mathcal P^{(3)}=\mathcal{SP}$.
The diagram of a $d$-dimensional partition $\alpha$ is the finite ideal
\[
 D_{d+1}(\alpha)=\{(i_1,\ldots,i_d,j):1\leq j\leq\alpha_{i_1,\ldots,i_d}\}
 \subseteq\ZZ_{>0}^{d+1}.
\]
Meet and join are componentwise minimum and maximum, volume is a valuation, and the nonzero join-irreducibles are the partitions whose diagrams are principal boxes.

For $r=d+1$, the sublattice generated by the $d$-dimensional partitions of volume at most $n$ is the ideal lattice of
\[
 Q_n^{(r)}=\{(i_1,\ldots,i_r)\in\ZZ_{>0}^r:i_1\cdots i_r\leq n\}.
\]
\editnov{The parallel box-bounded lattice is
\[
 Q_{\mathbf b}=[b_1]\times\cdots\times[b_r],\qquad
 L_{\mathbf b}^{(r)}=\Id(Q_{\mathbf b}),\qquad \mathbf b=(b_1,\ldots,b_r).
\]
We identify its ideals with $(r-1)$-dimensional partitions through their diagrams. Its toric algebra $A_{Q_{\mathbf b}}$ has dimension $b_1\cdots b_r$ and a basis indexed by chains of these partitions. The inclusions
\[
 Q_n^{(r)}\subseteq[n]^r\subseteq Q_{n^r}^{(r)},\qquad
 Q_{\mathbf b}\subseteq Q_{b_1\cdots b_r}^{(r)}
\]
show that the two families are cofinal among finite diagram regions. Extending an ideal by the same subset of the larger region gives compatible lattice and algebra inclusions, since each region is an order ideal. Consequently,
\[
 \varinjlim_n A_{Q_n^{(r)}}\simeq
 \varinjlim_{\mathbf b}A_{Q_{\mathbf b}},\qquad
 \varprojlim_n U_{Q_n^{(r)}}\simeq
 \varprojlim_{\mathbf b}U_{Q_{\mathbf b}},
\]
where $U_{Q_{\mathbf b}}=\Spec A_{Q_{\mathbf b}}$. The isomorphisms preserve volume and slice number.}

The maximum of the volume-generated lattice has entries $\lfloor n/(i_1\cdots i_d)\rfloor$. Indeed, any entry of a partition of volume at most $n$ satisfies
$(i_1\cdots i_d)\alpha_{i_1,\ldots,i_d}\leq n$, and every principal box in $Q_n^{(r)}$ has volume at most $n$. The same proof as Proposition~\ref{prop:solidfinite} applies.
Let $A_n^{(r)}=A_{Q_n^{(r)}}$ and $U_n^{(r)}=\Spec A_n^{(r)}$. These normal toric varieties have dimension
\[
 d_r(n)=|Q_n^{(r)}|=\sum_{m=1}^n\tau_r(m)\sim\frac{n(\log  n)^{r-1}}{(r-1)!},\qquad
 \tau_r(m)=\#\{(i_1,\ldots,i_r):i_1\cdots i_r=m\}.
\]
Their exponent cones consist of nonnegative order-reversing vectors on $Q_n^{(r)}$. Their dual cones are generated, as before, by $e_u-e_v$ for covers $u\lessdot v$ and by $e_u$ for maximal cells.

{

By induction on the dimension \(r\), the volume of
\[
\{x_1,\dots,x_k\ge1:x_1\cdots x_r\le n\}
\]
yields
\[
\sum_{\substack{i_1i_2\cdots i_r\le n\\i_1,\dots,i_r\ge1}}1
\sim\frac{n(\log  n)^{r-1}}{(r-1)!}.
\]
This is valid for \(r=2\) the asymptotic is \( n(\log  n)\).
}

\begin{proposition}[Dimension-independent extension]\label{prop:all-dimensions}
For every fixed $d\geq1$, the direct-limit algebra
$A_\infty^{(d+1)}=\varinjlim_n A_n^{(d+1)}$ has a standard monomial basis indexed by $(d+1)$-dimensional partitions, and
\[
 \Hilb(A_\infty^{(d+1)};q)=F_{d+1}(q).
\]
Moreover, lattice points in the $s$th dilation of the integral polytope
\[
 \rev{\Delta}_n^{(d+1)}=
 \operatorname{conv}\{\chi_I:\varnothing\ne I\in\Id(Q_n^{(d+1)})\}
\]
correspond to $(d+1)$-dimensional partitions with exactly $s$ nonzero $d$-dimensional slices in that finite lattice. The sum of the lattice coordinates is the volume. In each fixed volume $m$, the finite algebra stabilizes to the limit once $n\geq\max\{1,m\}$.
\end{proposition}
\begin{proof}
A chain of nonempty $d$-dimensional partitions is exactly a $(d+1)$-dimensional partition, sliced in one coordinate direction. Standard monomials are indexed by such chains. The sum of the corresponding ideal indicators gives the lattice-point bijection, with inverse given by superlevel ideals; the coordinate at the minimum cell records the number of slices. Finally, every slice of a volume-$m$ partition has volume at most $m$, which proves stabilization.
\end{proof}
For $n\geq2$, the standard Proj, its section spaces, and equivariant Riemann--Roch give the same geometric interpretation as in Section~\ref{sec:geometry}. Thus higher dimension introduces no additional hypothesis into these arguments.

\begin{example}[Hypercubes in every dimension]
For $r\geq2$, the poset $Q_2^{(r)}$ consists of its minimum and the $r$ cells obtained by changing one coordinate from one to two. Consequently
\[
 \rev{\Delta}_2^{(r)}=\{1\}\times[0,1]^r\subseteq\RR^{r+1},
 \qquad \#(s\rev{\Delta}_2^{(r)}\cap\ZZ^{r+1})=(s+1)^r.
\]
This counts $r$-dimensional partitions with exactly $s$ nonzero slices in the corresponding finite lattice. The standard Proj is the Segre variety $(\mathbb P^1)^r$, and the volume character in degree $s$ is $q^s(1+q+\cdots+q^s)^r$.
\end{example}

The signed identity also extends. Let $G_\infty^{(d+1)}$ be the incomparability graph of the nonempty elements of $\mathcal P^{(d)}$. Define $S_{d+1}$ by the same rule as $S_3$: retain the multisets along the slicing direction, require pairwise distinct slices, each weakly decreasing in its own coordinates, and arrange them in decreasing lexicographic order, and sum $(-1)^{|D|}$ over selections of incomparable slice pairs covering every slice. Set $S_{d+1}(\varnothing)=1$. Then
\begin{equation}\label{eq:all-dimensions}
 \sum_{\Pi\in\mathcal P^{(d+1)}}S_{d+1}(\Pi)q^{|\Pi|}
 =F_{d+1}(q)\prod_{m\geq1}(1-q^m)^{p_d(m)}.
\end{equation}
Indeed, the volume-compatible Gr\"obner degeneration gives the edge ideal of $G_\infty^{(d+1)}$, and the subgraph formula groups by canonically ordered slice sets exactly as above. Sorting fibers preserves all the remaining coordinate inequalities. Every bounded volume involves finitely many slices and graphs, so these steps are valid in all fixed dimensions. For $d=1$ this is the main theorem, with $S_2=S$; for $d=2$ it is \eqref{eq:solididentity}.

What does not extend is a closed evaluation of the ordinary-volume series by the classical products. Euler evaluates $F_1$, and MacMahon evaluates $F_2$. The analogous conjectural product in dimension $r$,
\[
 M_r(q)=\prod_{m\geq1}(1-q^m)^{-\binom{m+r-2}{r-1}},
\]
fails as an ordinary-volume formula for $r\geq3$; see \cite{Atkin,AY}. Formula~\eqref{eq:all-dimensions} remains exact without substituting $M_{d+1}$ for $F_{d+1}$. In addition, for $d\geq3$ its exponents $p_d(m)$ themselves count higher-dimensional partitions and lack the plane-partition product evaluation. The theory therefore extends uniformly, while explicit unrestricted enumeration remains the unresolved part. This limitation does not affect the finite polytope counts or their standard-Proj interpretation.

\subsection{The refined factorization in every dimension}\label{subsec:higher-refinement}
The refinement of Section~\ref{subsec:toric-refinement} also extends: its positive numerator comes from the exponent cone, independently of a MacMahon evaluation. We first state the finite result for any nonempty finite poset $Q$ with a least element $o$. Put $m=|Q|$, and let $\Gamma_Q$ be the incomparability graph on $\Id(Q)\setminus\{\varnothing\}$. Define
\[
 F_Q(t,q)=\sum_{h\in C_Q\cap\ZZ^Q}t^{h_o}q^{\ell(h)},
 \qquad \ell(h)=\sum_{u\in Q}h_u,
 \qquad p_Q(k)=\#\{I\in\Id(Q):|I|=k\}.
\]
%\textcolor{blue}{$h\in C_Q\cap \ZZ^{|Q|}$}

Under the monomial realization of $A_Q$, $h_o$ is the number of nonempty ideals in its standard chain and $\ell(h)$ is their total cardinality. For the graph numerator write
\[
 \mathcal N_Q(t,q)=
 \sum_{D\subseteq E(\Gamma_Q)}(-1)^{|D|}
 t^{|V(D)|}q^{\sum_{I\in V(D)}|I|},
\]
where $V(D)$ is the set of endpoints of the selected edges; the empty selection contributes $1$. The calligraphic notation distinguishes this polynomial from the dual lattice $N_Q$.
Choose a natural labeling $\nu:Q\to[m]$, meaning that $u<v$ implies $\nu(u)<\nu(v)$, and set
\[
 E_Q(t,q)=\sum_{w\in\operatorname{LinExt}(Q)}
 t^{\operatorname{des}(w)}q^{\operatorname{maj}(w)}.
\]
Here extensions are words of labels, with the descent statistics defined in Section~\ref{subsec:toric-refinement}.

\begin{proposition}[Refined factorization for finite toric approximations]\label{prop:higher-refined}
With the preceding notation,
\begin{align}
 F_Q(t,q)
 &=\frac{\mathcal N_Q(t,q)}
 {\prod_{\varnothing\ne I\in\Id(Q)}(1-tq^{|I|})}
 =\frac{E_Q(t,q)}{\prod_{k=1}^{m}(1-tq^k)},\label{eq:higher-refined-finite}\\
 \mathcal N_Q(t,q)
 &=E_Q(t,q)\prod_{k=1}^{m}(1-tq^k)^{p_Q(k)-1}.
 \label{eq:higher-refined-factor}
\end{align}
In particular the signed numerator factors through a polynomial with nonnegative coefficients. For $r\geq2$ and a finite order ideal $Q\subset\ZZ_{>0}^r$, the same numerator is
\[
 \mathcal N_Q(t,q)=\sum_{\Pi}S_r(\Pi)t^{\rho(\Pi)}q^{|\Pi|},
\]
where $\Pi$ ranges over the $r$-dimensional partitions whose nonempty $(r-1)$-dimensional slices have diagrams in $\Id(Q)$, and $\rho(\Pi)$ is their number. The slicing direction is fixed throughout.
\end{proposition}
\begin{proof}
The Hibi initial ideal is the edge ideal of $\Gamma_Q$. Since all nonempty ideals contain $o$, both standard degree and cardinality weight are preserved. The graph expansion, with vertex weight $tq^{|I|}$, gives the first quotient. When $Q$ is an ideal of a product grid, grouping edge selections by the sorted slices gives $S_r$ exactly as in \eqref{eq:all-dimensions}; it preserves both slice number and volume.

For the second quotient, order the coordinates of $h$ decreasingly and break ties by increasing label. This selects a unique extension $w$. Put $d_k=h_{w_k}-h_{w_{k+1}}$, with $h_{w_{m+1}}=0$. The conditions are $d_k\geq0$, with $d_k\geq1$ required at each descent. Because $o$ comes first,
\[
 h_o=\sum_k d_k,\qquad \ell(h)=\sum_k k d_k.
\]
The contribution is $t^{\operatorname{des}(w)}q^{\operatorname{maj}(w)}/\prod_k(1-tq^k)$. This is the classical linear-extension decomposition of $P$-partitions \cite[\S3.15, Lemma~3.15.3 and Theorem~3.15.5]{StanleyEC1}, or the half-open version of the order-polytope triangulation \cite[\S5]{StanleyPolytopes}. Summation and comparison of denominators prove the formulas. Every cardinality $1\leq k\leq m$ occurs among the ideals (take prefixes of an extension), so $p_Q(k)-1\geq0$.
\end{proof}

\begin{example}[A solid-partition approximation]\label{ex:solid-refined}
Take the four-cell poset
\[
Q=Q_2^{(3)}=\{(1,1,1),(2,1,1),(1,2,1),(1,1,2)\},
\]
and label these cells $1,2,3,4$, respectively. Every extension starts with $1$, followed by an arbitrary permutation of $2,3,4$. Hence
\[
 E_Q(t,q)=1+2tq^2+2tq^3+t^2q^5,
 \qquad (p_Q(1),p_Q(2),p_Q(3),p_Q(4))=(1,3,3,1),
\]
\[
 \mathcal N_Q(t,q)
 =(1+2tq^2+2tq^3+t^2q^5)(1-tq^2)^2(1-tq^3)^2.
\]
In coordinates the cone is $h_o=s\geq h_2,h_3,h_4\geq0$, so the same refined count is
\[
 F_Q(t,q)=\sum_{s\geq0}t^s q^s(1+q+\cdots+q^s)^3
 =\frac{1+2tq^2+2tq^3+t^2q^5}{\prod_{k=1}^{4}(1-tq^k)}.
\]
Thus the factorization gives an explicit solid-partition count at a finite toric stage without using a product formula for unrestricted solid partitions.
\end{example}

For the limit, use the cofinal boxes $Q_{\mathbf b}=[b_1]\times\cdots\times[b_r]$ \editnov{introduced in Section~\ref{subsec:all-dimensions}} and write $E_{r,\infty}=\lim_{\mathbf b\to\infty}E_{Q_{\mathbf b}}$. This limit exists coefficientwise in $\ZZ_{\geq0}[t][[q]]$: in volume at most $M$, all diagrams involved fit once every $b_i\geq M$, and \eqref{eq:higher-refined-finite} then gives stabilization of the numerator. If
\[
 F_r(t,q)=\sum_{\Pi\in\mathcal P^{(r)}}t^{\rho(\Pi)}q^{|\Pi|},
 \qquad
 \mathcal N_r(t,q)=\sum_{\Pi\in\mathcal P^{(r)}}S_r(\Pi)t^{\rho(\Pi)}q^{|\Pi|},
\]
then $F_r(1,q)=F_r(q)$ and
\begin{equation}\label{eq:higher-refined-limit}
 F_r(t,q)=\frac{E_{r,\infty}(t,q)}{\prod_{k\geq1}(1-tq^k)},
 \qquad
 \mathcal N_r(t,q)=E_{r,\infty}(t,q)
 \prod_{k\geq1}(1-tq^k)^{p_{r-1}(k)-1}.
\end{equation}
For $r=2$ these are Corollary~\ref{cor:refined-limit}, with $E_{2,\infty}=E_\infty$. For $r\geq3$ the factorization and positivity remain valid; only the explicit evaluation of $F_r(1,q)$, and hence of $E_{r,\infty}(1,q)$, is unavailable from the classical MacMahon product. The linear-extension formula is classical; its combination here with the signed slice statistic makes the recursive toric interpretation compatible with the refined graph numerator.


\begin{thebibliography}{99}
\bibitem{Afsharijoo} P. Afsharijoo, \emph{Even--odd partition identities of Rogers--Ramanujan type}, Ramanujan J. \textbf{57} (2022), 969--979. \url{https://doi.org/10.1007/s11139-021-00470-3}.

\bibitem{ADJM} P. Afsharijoo, J. Dousse, F. Jouhet, and H. Mourtada, \emph{New companions to the Andrews--Gordon identities motivated by commutative algebra}, Adv. Math. \textbf{417} (2023), 108946. \url{https://doi.org/10.1016/j.aim.2023.108946}.

\editnov{\bibitem{AGPM} P. Afsharijoo, P. D. Gonz\'alez P\'erez, and H. Mourtada, \emph{Partition identities associated with $A_r$ surface singularities}, preprint, arXiv:2601.12048, 2026. \url{https://arxiv.org/abs/2601.12048}.}

\bibitem{AM} P. Afsharijoo and H. Mourtada, \emph{Partition identities and application to infinite-dimensional Gr\"obner basis and vice versa}, in \emph{Arc Schemes and Singularities}, World Scientific, 2020, 145--161. \url{https://doi.org/10.1142/9781786347206_0009}.

\editnov{\bibitem{AMNeighborly} P. Afsharijoo and H. Mourtada, \emph{Neighborly partitions, hypergraphs and Gordon's identities}, J. Combin. Theory Ser. A \textbf{210} (2025), Paper No.~105963, 16 pp. \url{https://arxiv.org/abs/2309.13334}.}

\bibitem{AY} A. Amanov and D. Yeliussizov, \emph{MacMahon's statistics on higher-dimensional partitions}, Forum Math. Sigma \textbf{11} (2023), e63. \url{https://doi.org/10.1017/fms.2023.61}.

\bibitem{Andrews} G. E. Andrews, \emph{The Theory of Partitions}, Cambridge Mathematical Library, Cambridge University Press, 1998; reprint of the 1976 original.

\bibitem{Atkin} A. O. L. Atkin, P. Bratley, I. G. Macdonald, and J. K. S. McKay, \emph{Some computations for $m$-dimensional partitions}, Proc. Cambridge Philos. Soc. \textbf{63} (1967), 1097--1100. \url{https://doi.org/10.1017/S0305004100042171}.

\bibitem{BGS} R. Barman, A. Ghosh, and G. Singh, \emph{On the G\"ollnitz--Gordon--Andrews identities via commutative algebra}, Proc. Amer. Math. Soc., published electronically October 28, 2026. \url{https://doi.org/10.1090/proc/17932}.

\bibitem{Baxter} R. J. Baxter, \emph{Hard hexagons: Exact solution}, J. Phys. A \textbf{13} (1980), no.~3, L61--L70.

\bibitem{BR} M. Beck and S. Robins, \emph{Computing the Continuous Discretely: Integer-Point Enumeration in Polyhedra}, second edition, Undergraduate Texts in Mathematics, Springer, New York, 2015. \url{https://doi.org/10.1007/978-1-4939-2969-6}.

\bibitem{Brion} M. Brion, \emph{Points entiers dans les poly\`edres convexes}, Ann. Sci. \'Ecole Norm. Sup. (4) \textbf{21} (1988), no.~4, 653--663. \url{https://doi.org/10.24033/asens.1572}.

\bibitem{BMS} C. Bruschek, H. Mourtada, and J. Schepers, \emph{Arc spaces and the Rogers--Ramanujan identities}, Ramanujan J. \textbf{30} (2013), no.~1, 9--38. \url{https://doi.org/10.1007/s11139-012-9401-y}.

\bibitem{CLS} D. A. Cox, J. B. Little, and H. K. Schenck, \emph{Toric Varieties}, Graduate Studies in Mathematics 124, American Mathematical Society, 2011.

\bibitem{FF} B. Feigin and E. Frenkel, \emph{Coinvariants of nilpotent subalgebras of the Virasoro algebra and partition identities}, in \emph{I. M. Gel'fand Seminar}, Advances in Soviet Mathematics 16, American Mathematical Society, 1993, 139--148.

{\bibitem{FM} E. Feigin and I. Makhlin, \emph{Relative poset polytopes and semitoric degenerations}, Selecta Math. (N.S.) \textbf{30} (2024), Paper No.~48. \url{https://doi.org/10.1007/s00029-024-00935-5}.}

\bibitem{Fulman} J. Fulman, \emph{A probabilistic proof of the Rogers--Ramanujan identities}, Bull. London Math. Soc. \textbf{33} (2001), no.~4, 397--407. \url{https://doi.org/10.1017/S0024609301008207}.

\bibitem{GB} A. Ghosh and R. Barman, \emph{$J$-generalization of the Rogers--Ramanujan--Gordon identities via commutative algebra}, arXiv:2602.02187, 2026. \url{https://arxiv.org/abs/2602.02187}.

{\bibitem{GL96} N. Gonciulea and V. Lakshmibai, \emph{Degenerations of flag and Schubert varieties to toric varieties}, Transform. Groups \textbf{1} (1996), 215--248.}

{\bibitem{GL97} N. Gonciulea and V. Lakshmibai, \emph{Schubert varieties, toric varieties and ladder determinantal varieties}, Ann. Inst. Fourier (Grenoble) \textbf{47} (1997), no.~4, 1013--1064. \url{https://doi.org/10.5802/aif.1590}.}

\bibitem{Goodarzi} A. Goodarzi, \emph{On the Hilbert series of monomial ideals}, J. Combin. Theory Ser. A \textbf{120} (2013), 315--317.

\bibitem{GMMR} M. Graffeo, S. Monavari, R. Moschetti, and A. T. Ricolfi, \emph{Enumeration of partitions via socle reduction}, preprint, arXiv:2501.10267v1, 2025. \url{https://arxiv.org/abs/2501.10267v1}.

\bibitem{HH} J. Herzog and T. Hibi, \emph{Monomial Ideals}, Graduate Texts in Mathematics 260, Springer, 2011. \url{https://doi.org/10.1007/978-0-85729-106-6}.

\bibitem{HHO} J. Herzog, T. Hibi, and H. Ohsugi, \emph{Binomial Ideals}, Graduate Texts in Mathematics 279, Springer, 2018.

\bibitem{Hibi} T. Hibi, \emph{Distributive lattices, affine semigroup rings and algebras with straightening laws}, in \emph{Commutative Algebra and Combinatorics}, M. Nagata and H. Matsumura (eds.), Advanced Studies in Pure Mathematics, Vol.~11, North-Holland, Amsterdam, 1987, 93--109.

\bibitem{LM} V. Lakshmibai and H. Mukherjee, \emph{Singular loci of Hibi toric varieties}, J. Ramanujan Math. Soc. \textbf{26} (2011), no.~1, 1--29. \url{https://arxiv.org/abs/0706.2686}.

\bibitem{MacMahon1896} P. A. MacMahon, \emph{Memoir on the theory of the partition of numbers---Part I}, Philos. Trans. Roy. Soc. London Ser. A \textbf{187} (1896), 619--673. \url{https://doi.org/10.1098/rsta.1896.0016}.

\bibitem{MacMahon} P. A. MacMahon, \emph{Combinatory Analysis}, Vol.~II, Cambridge University Press, 1916.

\bibitem{MS} E. Miller and B. Sturmfels, \emph{Combinatorial Commutative Algebra}, Graduate Texts in Mathematics 227, Springer, 2005.

\bibitem{MM} Z. Mohsen and H. Mourtada, \emph{Neighborly partitions and the numerators of Rogers--Ramanujan identities}, Int. J. Number Theory \textbf{19} (2023), no.~4, 859--872. \url{https://doi.org/10.1142/S1793042123500434}.

\bibitem{Mourtada} H. Mourtada, \emph{Hilbert meets Ramanujan: Singularity theory and integer partitions}, Bull. Amer. Math. Soc. (N.S.) \textbf{62} (2025), no.~1, 93--111. \url{https://doi.org/10.1090/bull/1854}.

\bibitem{OR} A. Okounkov and N. Reshetikhin, \emph{Correlation function of Schur process with application to local geometry of a random 3-dimensional Young diagram}, J. Amer. Math. Soc. \textbf{16} (2003), 581--603. \url{https://doi.org/10.1090/S0894-0347-03-00425-9}.

\bibitem{ORV} A. Okounkov, N. Reshetikhin, and C. Vafa, \emph{Quantum Calabi--Yau and classical crystals}, in \emph{The Unity of Mathematics}, P. Etingof, V. Retakh, and I. M. Singer (eds.), Progress in Mathematics, Vol.~244, Birkh\"auser, Boston, 2006, 597--618. \url{https://arxiv.org/abs/hep-th/0309208}.

\bibitem{Peeva} I. Peeva, \emph{Graded Syzygies}, Algebra and Applications, Vol.~14, Springer-Verlag, London, 2011, xii+302 pp.

\bibitem{SA} F. Schreier-Aigner, \emph{Fully complementary higher dimensional partitions}, Ann. Comb. \textbf{29} (2025), 1--23. \url{https://doi.org/10.1007/s00026-024-00691-5}.

\bibitem{StanleyPP} R. P. Stanley, \emph{Theory and application of plane partitions}, Parts 1 and 2, Stud. Appl. Math. \textbf{50} (1971), 167--188 and 259--279.

\bibitem{StanleyHilbert} R. P. Stanley, \emph{Hilbert functions of graded algebras}, Adv. Math. \textbf{28} (1978), no.~1, 57--83. \url{https://doi.org/10.1016/0001-8708(78)90045-2}.

\bibitem{StanleyPolytopes} R. P. Stanley, \emph{Two poset polytopes}, Discrete Comput. Geom. \textbf{1} (1986), 9--23.

\bibitem{StanleyEC1} R. P. Stanley, \emph{Enumerative Combinatorics}, Vol.~1, second edition, Cambridge Studies in Advanced Mathematics 49, Cambridge University Press, 2012, \S\S1.4 and 3.15. \url{https://doi.org/10.1017/CBO9781139058520}.

\bibitem{StanleyEC2} R. P. Stanley, \emph{Enumerative Combinatorics}, Vol.~2, Cambridge Studies in Advanced Mathematics 62, Cambridge University Press, 1999, Chapter~7.
\bibitem{Vuletic} M. Vuleti\'c, \emph{A generalization of MacMahon's formula}, Trans. Amer. Math. Soc. \textbf{361} (2009), 2789--2804.
\end{thebibliography}
\end{document}